# SMARANDACHE SPECIAL DEFINITE ALGEBRAIC STRUCTURES

**W. B. Vasantha Kandasamy**
e-mail: **vasanthakandasamy@gmail.com**
web: **http://mat.iitm.ac.in/~wbv**
**www.vasantha.in**

**2009**

# SMARANDACHE SPECIAL DEFINITE ALGEBRAIC STRUCTURES

W. B. Vasantha Kandasamy

**2009**

# CONTENTS

# PREFACE

In this book we introduce the notion of Smarandache special definite algebraic structures. We can also call them equivalently as Smarandache definite special algebraic structures. These new structures are defined as those strong algebraic structures which have in them a proper subset which is a weak algebraic structure. For instance, the existence of a semigroup in a group or a semifield in a field or a semiring in a ring. It is interesting to note that these concepts cannot be defined when the algebraic structure has finite cardinality i.e., when the algebraic structure has finite number of elements in it.

This book has four chapters. Chapter one is introductory in nature. In chapter two the notion of Smarandache special definite groups and Smarandache special definite fields are introduced and several interesting properties are derived. The notion of Smarandache definite special rings, vector spaces and linear algebras are introduced and analysed in chapter three. The final chapter suggests over 200 problems.

I deeply acknowledge the unflinching support of Dr.K.Kandasamy, Meena and Kama.

W.B.VASANTHA KANDASAMY

# DEDICATION

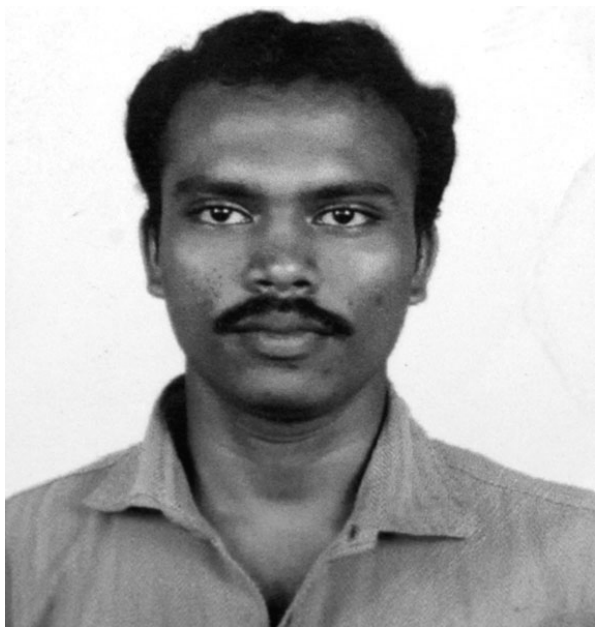

**Muthukumar Kumaresan**
19.11.1982 – 29.01.2009

*We dedicate this book to the 'Heroic Tamil Son' Muthukumar Kumaresan who immolated himself in Chennai, Tamil Nadu on 29 January 2009 to draw the world's attention to the genocide of Tamils in Sri Lanka. Muthukumar's supreme self-sacrifice brought about an unparalleled students' upsurge in
Tamil Nadu, and became the rallying point for the struggle to recognize Tamil Eelam, a separate homeland for Tamils in the embattled island of Sri Lanka.
We salute the courage of this braveheart.*

**Chapter One**

# BASIC CONCEPTS

In this chapter we just recall some of the basic notions essential to make this book a self contained one. This chapter has three sections. In the first section notions about groups and semigroups are introduced. Section two introduces the concepts of fields and rings. Section three introduces the notion of semirings and semifields. The final section recalls the concepts of vector spaces and semivector spaces.

In this chapter we recall the notion of groups and semigroups. We know all groups are semigroups and semigroups in general are not groups. That is why one can always think of semigroups as a generalization of groups. We have studied Smarandache semigroups [100], which are a new class of algebraic structures which contain a proper subset P of S such that P is a group. Thus Smarandache semigroups are semigroups which has a proper subset P of S such that P is a group under the operations of the semigroup.

Now in this book we study a new algebraic Smarandache structure called Smarandache definite algebraic structures which contains a proper subset which satisfies an algebraic structure under the same operations but is weaker than the given algebraic structure i.e., for instance finding proper subsets which are semigroups in the group. Like wise for a field to contain proper subsets which are rings, and rings which contain subsets which are semirings and so on. Finally we study those structures which are vector spaces that contain modules as substructures.

## 1.1 Groups, Semigroups and Smarandache Semigroups

In this section we just recall the basic definition of groups, semigroups and Smarandache semigroups. For more about these concepts the interested reader can refer [100].

**DEFINITION 1.1.1:** *Let $(G, *)$ denote a non empty set $G$ under the associative binary operation $*$. We call $(G, *)$ to be a group if the following conditions are satisfied.*

1. *$G$ is closed under the binary operation $*$ i.e., for all $a, b \in G$, $a * b \in G$.*
2. *There exists a unique element $e$ in $G$ such that $a * e = e * a = a$ for all $a \in G$ called the identity element of $G$.*
3. *There exists for every $a \in G$ a unique element $a' \in G$ such that $a * a' = a' * a = e$; $a'$ is defined as the inverse of $a$ in $G$.*

If '$*$' is a binary operation on $G$ such that $a * b = b * a$ for all $a, b \in G$ we call $G$ to be a commutative group.

The number of distinct elements in $G$ is called the order of $G$ and denoted by $o(G)$ or $|G|$. If $o(G) = n$ and $n < \infty$ we call $G$ to be a finite group. If $o(G) = \infty$ then we call $G$ to be an infinite group.

If the group $G$ is generated by a single element $g$ then $G$ is said to be a cyclic group and is denoted as $G = \langle g \mid g^n = identity \rangle$.

We now illustrate these situations by the following examples.

*Example 1.1.1:* Let Z be the set of positive and negative integers with zero. Z is a group under addition. Z is a commutative group and order of Z is infinite.

*Example 1.1.2:* Let $Q^* = Q \setminus \{0\}$ be the set of all rationals. $Q^* = Q \setminus \{0\}$ is a group under multiplication. $Q^*$ is an infinite abelian group of infinite order.

*Example 1.1.3:* Let

$$M_{3\times 3} = \left\{ \begin{pmatrix} a & b & c \\ d & e & f \\ g & h & i \end{pmatrix} = A \mid a, b, c, ..., i \in Q \right\}$$

with $|A| \neq 0$. $M_{3\times 3}$ is a group under matrix multiplication. In fact $M_{3\times 3}$ is a non commutative infinite group.

*Example 1.1.4:* Let $S_n$ be the symmetric group of order n (n < ∞). $S_n$ be the symmetric group of order $\angle n = n!$.

*Example 1.1.5:* Let $G = \langle g \mid g^{16} = 1 \rangle$ be a group of order 16. G is a finite cyclic group.

*Example 1.1.6:* Let $Z_{26} = \{\bar{0}, \bar{1}, \bar{2}, ..., \overline{25}\}$, $Z_{26}$ is an additive group of order 26 which is finite and commutative.

*Example 1.1.7:* Let $D_{26} = \{a, b \mid a^2 = b^6 = 1, bab = a\}$, $D_{26}$ is the non commutative group of order 12 known as the dihedral group.

Now we request the reader to refer [35, 37, 100] for more about groups. Now we proceed on to recall the definition of semigroups.

**DEFINITION 1.1.2:** *Let (S, \*) be a non empty set on which is defined an associative binary operation \*. If for all a,b $\in$ S we have a \* b $\in$ S then we call S to be a semigroup with respect to the binary operation \* . If the number of distinct elements in S is finite then we call S to be a finite semigroup. If S has infinite number of elements then we call S to be an infinite semigroup. If* *in S, for all a, b in S and we have a * b = b * a then we call S to be a commutative semigroup. If a semigroup is generated by a single element then we call the semigroup to be a cyclic semigroup. A semigroup (S, *) in which there exists an element e $\in$ S with e * s = s * e = s for all s $\in$ S, we call e the identity element of S. S is called a monoid. Thus a semigroup with identity is called the monoid.*

Now we illustrate these situations by some examples.

*Example 1.1.8:* Let $Z^+$ be the set of positive integers. $Z^+$ is a semigroup under multiplication.

*Example 1.1.9:* Let Q be the set of rationals. Q is a monoid under multiplication for 1 $\in$ Q is such that 1 · q = q · 1 = q for all q $\in$ Q.

*Example 1.1.10:* Let $M_{2\times 2}$ = {$(a_{ij})$ = A such that $a_{ij} \in$ Q}; $M_{2\times 2}$ is a semigroup under matrix multiplication. In fact $M_{2\times 2}$ is an infinite monoid which is non commutative.

*Example 1.1.11:* Let $Z_{15}$ = { $\bar{0}, \bar{1}, \bar{2},\ldots,\overline{14}$ } be the set of integers modulo 15. $Z_{15}$ is a semigroup under multiplication. $Z_{15}$ is a commutative finite semigroup.

*Example 1.1.12:* Let $Z_{11}$ be the set of integers modulo 11. $Z_{11}$ is a finite commutative semigroup under addition modulo 11. $Z_{11}$ is in fact a cyclic semigroup for it is generated by 1.

Now we have seen the examples; we just recall the definition of subsemigroup and ideal.

**DEFINITION 1.1.3:** *Let (S, *) be a semigroup. A proper non empty subset P of S is said to be a subsemigroup of S if (P, *) is a semigroup.*

We illustrate this by the following example.

*Example 1.1.13:* Let $Z^+$ be the semigroup of integers under multiplication. $3Z^+$ is a subsemigroup of $Z^+$ under multiplication. In fact $nZ^+$ for any $n \in Z^+$ is a subsemigroup of Z.

**DEFINITION 1.1.4:** *Let (S, \*) be a semigroup. Let I be a proper subset of S. I is called an ideal of S if the following conditions are satisfied.*
  1. *(I, \*) is a subsemigroup of (S, \*).*
  2. *For each $i \in I$ and $s \in S$; $i * s = s * i \in I$.*

We illustrate this by the following example:

*Example 1.1.14:* Let Z be the semigroup under multiplication. $Z^+$ is a subsemigroup of Z. Clearly $Z^+$ is not an ideal of Z for we see if $-5 \in Z$ and $3 \in Z^+$; $-5 \times 3 = -15$ is not in $Z^+$.

*Example 1.1.15:* Let $Z^+$ be the semigroup under multiplication. $13Z^+$ is a subsemigroup of $Z^+$ and in fact $13Z^+$ is an ideal of $Z^+$ under multiplication.

Now we have just seen from the example 1.1.14 that subsemigroups of a semigroup in general need not be an ideal of the semigroup.

Now we just recall the definition of a subgroup and the normal subgroup of a group.

**DEFINITION 1.1.5:** *Let (G, \*) be a group. A proper subset H of G is said to be a subgroup of G if (H, \*) is itself a group.*
  *We call a subgroup H of a group G to be a normal subgroup of G if $g * H * g^{-1} = H$ for all $g \in G$.*

We illustrate the above definition by the following example.

*Example 1.1.16:* Let $S_7$ be a group; the group of permutations of seven elements. $A_7$ the alternating subgroup of $S_7$. In fact $A_7$ is a normal subgroup of $S_7$.

*Example 1.1.17:* Let $Q \setminus \{0\}$ be the group under multiplication. $Q^+ \setminus \{0\}$ is a subgroup of $Q$ under multiplication. In fact $Q^+ \setminus \{0\}$ is also a normal subgroup of $Q \setminus \{0\}$.

All subgroups of a group G in general need not be normal subgroups of G.

We illustrate this by the following example.

*Example 1.1.18:* Let $S_4$ be the symmetric group of degree 4. Take

$$P = \left\{ \begin{pmatrix} 1 & 2 & 3 & 4 \\ 1 & 2 & 3 & 4 \end{pmatrix}, \begin{pmatrix} 1 & 2 & 3 & 4 \\ 2 & 3 & 4 & 1 \end{pmatrix}, \begin{pmatrix} 1 & 2 & 3 & 4 \\ 3 & 4 & 1 & 2 \end{pmatrix}, \begin{pmatrix} 1 & 2 & 3 & 4 \\ 4 & 1 & 2 & 3 \end{pmatrix} \right\}$$

a proper subset of $S_4$; P is a subgroup of $S_4$ ; P is not a normal subgroup of $S_4$.

Now we just recall the definition of Smarandache semigroup and give some examples.

**DEFINITION 1.1.6:** *Let $(S, \circ)$ be a semigroup we call S to be a Smarandache semigroup (S-semigroup) if S contains a non empty proper subset P such that $(P, \circ)$ is a group.*

We give an example of it.

*Example 1.1.19:* Let $Z_{10}$ be the semigroup under multiplication modulo 10. $\langle 3 \rangle = \{3, 9, 7, 1\} = P$

| × | 3 | 9 | 7 | 1 |
|---|---|---|---|---|
| 3 | 9 | 7 | 1 | 3 |
| 9 | 7 | 1 | 3 | 9 |
| 7 | 1 | 3 | 9 | 7 |
| 1 | 3 | 9 | 7 | 1 |

Clearly P is a group under multiplication modulo 10 is given in the above table. Hence $Z_{10}$ is a S-semigroup.

***Remark:*** If S is a S-semigroup such that S has only finite number of distinct elements then we call S to be a finite S semigroup. If S has infinite number of distinct elements then we call S to be a infinite S-semigroup. Let S be a semigroup, if S has a non empty subset P such that P has proper subset. T such that T is a group then we call P a Smarandache subsemigroup. It is interesting to note that any semigroup which has a S sub semigroup itself becomes a S-semigroup. It is still important to see that not all semigroups are S-semigroups. For take $Z^+$ a semigroup under multiplication. Clearly $Z^+$ has no proper subset which is a group under multiplication. Let Z be the semigroup under multiplication $S = \{-1, 1\}$ is a group under multiplication. Hence Z is a S-semigroup.

For several properties about this structure please refer [100].

### 1.2. Fields and Rings

Now we just recall the definition of fields and rings.

**DEFINITION 1.2.1:** *Let R be a set with two associative binary operations +, ×. We call R to be a ring if the following conditions are satisfied.*

1. *(R, +) is a group.*
2. *(R, ×) is a semigroup.*
3. *(i) $a \times (b + c) = a \times b + a \times c$*
   *(ii) $(a + b) \times c = a \times c + b \times c$*

*for all $a, b, c \in R$.*

*If R contains an element 1 such that $a \times 1 = 1 \times a = a$ for all $a \in R$ we call R to be a ring with unit. If in R, $a \times b = b \times a$ for all $a, b \in R$ then we call R to be a commutative ring. If in the commutative ring R we have $a \times b = 0$ implies either $a = 0$ or $b = 0$ then we call R to be an integral domain. If $(R, +, \times)$ be an integral domain and if in addition $(R \setminus \{0\}, \times)$ is a commutative group then we call $(R, +, \times)$ to be a field. Clearly all fields are*

*commutative rings with unit. But all commutative rings with unit need not be a field. Suppose (R, +, ×) be a ring with unit and if (R \ {0}, ×) is a non commutative group. Then we call (R, +, ×) to be a division ring. All division rings are rings but a ring in general is not a division ring.*

Now we will proceed on to illustrate these by examples.

***Example 1.2.1:*** Let $M_{2\times 2} = \{M = (a_{ij}) \mid a_{ij} \in Q\}$ be the collection of all $2 \times 2$ matrices with entries from Q. $M_{2\times 2}$ is a ring which is non commutative but $M_{2\times 2}$ is not a division ring. Thus every non commutative ring need not be a division ring.
    Clearly all division rings are rings.

***Example 1.2.2:*** Let $Z_{12}$ be the ring of integers module 12. $Z_{12}$ is a finite commutative ring.

***Example 1.2.3:*** Let $M_{2\times 2}^2 = \{(a_{ij}) = M \mid a_{ij} \in \{0, 1\}\}$. $M_{2\times 2}^2$ denotes the set of all $2 \times 2$ matrices with entries from the prime field of characteristic two. $M_{2\times 2}^2$ is a finite non commutative ring. In fact it is not even a division ring.

***Example 1.2.4:*** Let $Z_5 = \{\bar{0}, \bar{1}, \bar{2}, \bar{3}, \bar{4}\}$ be the set of integers modulo 5. $Z_5$ is a field. It is in fact a finite field.

***Example 1.2.5:*** Let Z be the ring of integers. Z is a commutative ring, in fact an integral domain. 3Z is also an integral domain which is not a field.

***Example 1.2.6:*** Let Q be the ring of rationals. Q is a field. R the set of reals is also a field. C the collection of complex numbers is also a field. We see $Q \subseteq R \subseteq C$ even as sets.

***Example 1.2.7:*** Let Q[x] be the ring of polynomials. Q[x] is only a ring and not a field known as the ring of polynomials.

***Example 1.2.8:*** Let $Z_7[x]$ be the ring of polynomials in the variable x.

$$Z_7[x] = \left\{ \sum_{i=0}^{n} \alpha_i x_i \,\middle|\, n \in N; \alpha_i \in Z_7 \right\}.$$ Clearly $Z_7[x]$ is not a field only a commutative ring with unit.

***Example 1.2.9:*** Let $Z_6[x]$ be the ring of polynomials in the indeterminate x, $Z_6[x]$ is a commutative ring having zero divisiors.

Now we proceed on to define the notion of subrings and ideals of a ring.

**DEFINITION 1.2.2:** *Let R be any ring. A proper non empty subset P of R is said to be a subring of R if P is a ring under the operations of R.*

**DEFINITION 1.2.3:** *Let R be any ring. A proper non empty subset I of R is said to be a ideal of R if the following conditions are satisfied.*
  1. *I is a subring of R.*
  2. *For every $i \in I$ and $r \in R$; ir and ri $\in I$.*

Thus from the very definitions we see all ideals are subrings and subrings in general need not be ideals of R.

***Example 1.2.10:*** Let Q be the ring. Z is a subring of Q but Z is not an ideal of Q.

***Example 1.2.11:*** Let $Z_{12}[x]$ be the polynomial ring. $Z_{12}$ is a subring of $Z_{12}[x]$. Clearly $Z_{12}$ is not an ideal of $Z_{12}[x]$.

***Example 1.2.12:*** Let Z be the ring of integers. 5Z is a subring of Z as well as ideal of Z.

***Example 1.2.13:*** Let $M_{2 \times 2} = \{(m_{ij}) = M \mid m_{ij} \in Q\}$ be the ring of matrices. $P_{2 \times 2} = \{(m_{ij}) = M \mid m_{ij} \in Z\}$ is a subring of $M_{2 \times 2}$ which is not an ideal of $M_{2 \times 2}$.

Having seen the examples of subrings and ideals in a ring we now proceed on to make a mention of special types of ideals for more about these concepts please refer [48, 60].

We have defined an ideal. Let R be any ring. I ≠ (0) be an ideal of R. We say I is maximal ideal of R if we have any other ideal J of R such that (0) ≠ I ⊆ J ⊆ R, then either J = I or J = R. We call an ideal (0) ≠ P or R to be a minimal ideal of R if we have for any N an ideal of R such that (0) ⊆ N ⊆ P, either N = (0) or N = P. An ideal I of R is said to be prime if xy ∈ I then either x or y is in I. An ideal J of R is said to be principal if J is generated by a single element.

We illustrate these with examples.

*Example 1.2.14:* Let Z be the ring of integers 3Z is an ideal of Z. In fact 3Z is a maximal ideal of Z. 3Z is also a principal ideal of Z as 3Z is generated by 3. In fact 3Z is a prime ideal of Z as for xy ∈ 3Z clearly x or y is in 3Z. We see all ideals of the form pZ, p a positive prime is an ideal of Z which is maximal, prime and principal. We consider 6Z the ideal of Z, clearly 6Z is not a maximal ideal of Z · 6Z is not a prime ideal as 2 · 3 = 6 ∈ 6Z but both 2 and 3 are not in 6Z. 6Z is a principal ideal as it is generated by 6. We see
   2Z = {0, ±2, ±4, ±6, …}
   6Z = {0, ±6, ±12, ±18, … }
6Z ⊆ 2Z so 6Z is not a maximal ideal of Z. It is interesting to see Z has no minimal ideal.

*Example 1.2.15:* Let $Z_{12}$ be the ring of integers modulo 12. I = {0, 2, 4, 6, 8, 10}, J = {0, 6}, T = {0, 4, 8} and K = {0, 3, 6, 9} are ideals of $Z_{12}$. T and J are minimal ideals of $Z_{12}$, whereas I and K are maximal ideals of $Z_{12}$. In fact all the four ideals are principal. J and T are not prime ideals for 2 · 2 = 4 ∈ T but 2 ∉ T. Also 3·2 ∈ J but 2 and 3 ∉ J. Hence the claim.

Now we proceed on to give yet another example.

*Example 1.2.16:* Let $Z_6$ be the ring of integers modulo 6. I = {0, 3} and J = {0, 2, 4} are ideals of $Z_6$. I and J are both prime ideals of $Z_6$, they are also principal ideals of $Z_6$. It is very interesting to note that both I and J are maximal as well as minimal ideals of $Z_6$.

It is interesting to know a field cannot have proper ideals [35-8]. A field may have only subfields.

**DEFINITION 1.2.4:** *Let F be any field. We say F is a field of characteristic zero if nx = 0 for all x ∈ F; n an integer implies n = 0 i.e., for no integer nx = 0 for all x ∈ F can occur then F is a field of characteristic zero.*

*We say a field F is of characteristic p if we have px = 0 for all x ∈ F and p is a prime number.*

*Thus we say the field F is of characteristic p, p > 0 and p a prime.*

We illustrate this by the following examples.

*Example 1.2.17:* Let $Z_p = \{\bar{0}, \bar{1}, \bar{2}, \ldots, \overline{p-1}\}$ be the ring of integers modulo p. $Z_p$ is a field we see px = 0 for all x ∈ $Z_p$ so $Z_p$ is a field characteristic p.
   $Z_5$ is the field of characteristic 5. $Z_{13}$ is the field of characteristic 13. $Z_{12}$ is not a field but only a ring.

*Example 1.2.18:* Let R be the field of reals. R is of characteristic 0. Further every subset of R is also of characteristic zero. Q be the rationals; Q is the field of characteristic zero.

**DEFINITION 1.2.5:** *Let F be a field. A proper subset P of F is said to be a subfield of F if P is itself a field under the operations of F. If a field has no proper subset which is a subfield then we call F to be a prime field.*

*We see Q, the field of rationals is a prime field and Q is of characteristic zero. $Z_p$, p a prime are all prime fields of characteristic p. R the field of reals is not a prime field for it contains Q to be a subfield. Likewise C the field of complex numbers is not a prime field for it contains Q and R to be proper subfields.*

Now we define the notion of quotient rings for that alone leads us to show the existence of non prime fields of characteristic p. For we see all fields $Z_p$, p a prime are only prime fields of characteristic p. Let $Z_p[x]$ be the ring of polynomials with coefficients from $Z_p$. Let p(x) be an irreducible polynomial of degree n over $Z_p[x]$.

Consider the ideal I generated by p(x), $\frac{Z_p[x]}{I}$ is the quotient ring and if I is a maximal ideal the quotient ring $\frac{Z_p[x]}{I}$ is a field and $\frac{Z_p[x]}{I}$ has $p^n$ elements in it and $\frac{Z_p[x]}{I}$ is also a field of characteristic p and it has a subfield isomorphic to $Z_p$ given by $\{I, \bar{1}+I, \bar{2}+I, \ldots, \overline{p-1}\}$ where I acts as the additive zero of $\frac{Z_p[x]}{I}$.

We illustrate this situation by the following examples.

***Example 1.2.19:*** Let $Z_3[x]$ be the polynomial ring with coefficients in $Z_3$ in the variable x. Consider the polynomial $p(x) = 1 + x^4$ in $Z_3[x]$. Clearly p(x) is irreducible over $Z_3[x]$. Let I be the ideal generated by p(x).

It is easily verified $\frac{Z_3[x]}{I}$ has 81 elements in it and it is a field of characteristic 3. I acts as the additive identity and $\bar{1}+I$ acts as the multiplicative identity. For instance the multiplicative inverse of $2x^2 + I$ is for $(2x^2 + I)(x^2 + I) = \bar{1} + I$.

We now yet illustrate a few more examples.

***Example 1.2.20:*** Let $Z_2[x]$ be a polynomial ring with coefficients from $\{0, 1\}$ in the variable x. Consider a polynomial $p(x) = x^3 + x + 1$. Clearly p(x) is irreducible over $Z_2[x]$. Let I be the ideal generated by p(x). Consider the quotient ring

$$\frac{Z_2[x]}{\langle p(x)\rangle} = \frac{Z_2[x]}{I} = T$$
$$= \{I,\ \bar{1}+I,\ x+I,\ x^2+I,\ \bar{1}+x+I,\ \bar{1}+x^2+I,$$
$$x+x^2+I,\ \bar{1}+x+x^2+I\}.$$

Clearly T is a field of characteristic two. We have $S = \{I,\ \bar{1} + I\}$ is a proper subset of T which is a subfield of T. It is easily verified $S \cong Z_2 = \{0, 1\}$ for map $I \mapsto 0$ and $\bar{1}+I \mapsto 1$.

Now it is important to note that in general all quotient rings are not fields.

We illustrate this by the following examples.

*Example 1.2.21:* Let $Z_2[x]$ be the ring of polynomial over the prime field of characteristic two. Consider $p(x) = x^2 + 1 \in Z_2(x)$. Let J be the ideal generated by $x^2 + 1$. The quotient ring
$$\frac{Z_2[x]}{J} = \{J,\ \bar{1}+J,\ x+J,\ x+\bar{1}+J\} = P(\text{say}).$$
We prove P is not a field. Consider

$$\left((x+\bar{1})+J\right)\left((x+\bar{1})+J\right) = (x+\bar{1})^2 + J$$
$$= \left(x^2+2x+\bar{1}+J\right) = \left(x^2+\bar{1}+J\right) = J\ (\because x^2+\bar{1}\in J).$$

So $(x + \bar{1} + J)(x + \bar{1} + J) = J$. J is the zero of $\frac{Z_2[x]}{J} = P$ so $x + \bar{1} + J$ is a zero divisor in P. Hence P is a not a field it is only a quotient ring of characteristic two.

Now we have seen ideals and quotient rings. We just make a mention of how we have quotient rings to be also prime fields.

*Example 1.2.22:* Let Z be the ring of integers. Z is a ring, I = $\langle 3Z \rangle$ is an ideal of Z in fact a maximal ideal of Z. The quotient ring $\frac{Z}{I} = \{I, \bar{1}+I, \bar{2}+I\} = S$. Clearly S is isomorphic to the prime field of characteristic 3. For $Z_3 = \{\bar{0}, \bar{1}, \bar{2}\}$ then $S \cong Z_3$ given by $\phi(I) \mapsto \bar{0}$, $\phi(\bar{1}+I) \mapsto \bar{1}$ and $\phi(\bar{2}+I) \mapsto \bar{2}$. It is easily verified $\phi$ is a field isomorphism.

Now we give an example of a quotient ring which is not a prime field. We can generalise the example 1.2.21 as follows. Let Z be the ring of integers p be any prime. $pZ = I$ be the ideal generated by p. The quotient ring

$$\frac{Z}{pZ} = \frac{Z}{I} = V = \{I, \bar{1}+I, \bar{2}+I, \ldots, \overline{p-1}+I\}$$

and $V \cong Z_p$ ($Z_p$ the prime field of characteristic p).

Consider the following example which gives rings as quotient rings.

***Example 1.2.23:*** Let Z be the ring of integers. 12Z is an ideal generated by 12. Let I = 12Z, consider the quotient ring

$$\frac{Z}{I} = \{I, \bar{1}+I, \bar{2}+I, \ldots, \overline{11}+I\}.$$

We see $(\bar{3} + I)(\bar{4} + I) = \overline{12} + I = I$ is a zero divisor in the quotient ring $\frac{Z}{I} = T$. Also

$$(\bar{2} + I)(\bar{6} + I) = \overline{12} + I = I, (\bar{3} + I)(\bar{8} + I) = I.$$

Thus T is only a ring which is isomorphic to $Z_{12}$.

Now having given the recollection of definitions of rings and fields we now proceed on to give the definition of Smarandache rings, (S-rings).

**DEFINITION 1.2.6:** *Let R be a ring. We call R to be a Smarandache ring (S-ring) if we can find a proper subset P in R such that $P \neq \phi$ and $P \neq R$ but P is a field under the operations of R. All Smarandache rings are rings but in general every ring need not be a S-ring.*

We illustrate this situation by some examples.

*Example 1.2.24:* Let Q[x] be the polynomial ring over Q. We see Q the field of rationals is a proper subset of Q[x]. Thus Q[x] is a S - ring.

*Example 1.2.25:* Consider the polynomial ring Z[x]. Z[x] has no proper subset P which is a field. Hence Z[x] is not a S-ring. Thus we see all rings are not S-rings, but clearly every S-ring is a ring.

Now we proceed on to define the notion of Smarandache subring of any ring R.

**DEFINITION 1.2.7:** *Let R be any ring. P a proper subset of R. If P is a S-ring then we call P to be Smarandache subring (S-subring) of R.*

It is interesting to note that if a ring R has a proper Smarandache subring P then R itself is a S-ring; for P ⊆ R and P has a proper subset T such that T is a field now T ⊂ P ⊂ R so R itself is a S-ring.

*Example 1.2.26:* Let R[x] be the ring of polynomials over the reals R. Clearly Q[x] is the subring of R[x]. Q[x] is not a field but Q[x] is a S-ring for it contains a proper subset Q which is a field. Now Q ⊆ Q[x] ⊆ R[x] so R[x] is also a S-ring.

Now we proceed on to give example of a finite S-ring.

*Example 1.2.27:* Let $Z_2[x]$ be the ring of polynomials. Consider the quotient ring

$$\frac{Z_2[x]}{\langle x^2+1 \rangle = I} = \left\{ I,\ \bar{1}+I,\ x+I,\ \bar{1}+x+I \right\} = T.$$

T is a S-ring but T does not contain a S-subring. $Z_2[x]$ is also a S-ring. So every S-ring need not in general contain a S-subring. Now we recall the definition of Smarandache ideal of a ring.

**DEFINITION 1.2.8 [101]:** *The Smarandache ideal (S-ideal) is defined as an ideal A of a ring R such that a proper subset of A is a field(with respect to the same induced operations).*

It is interesting to note with these conditions we may not find examples of S-ideals. For $Z_6$ has no S-ideals, $Z_{10}$ has no S-ideals. So we proceed on to define Smarandache definite ideal.

**DEFINITION 1.2.9:** *The $(A, +, \cdot)$ be a S-ring. Let B be a proper subset of A which is a field. A non empty subset S of A is said to be Smarandache definite right ideal (S-definite right ideal) of A related to B if*
1. *$(S, +)$ is an additive abelian group.*
2. *For $b \in B$ and $s \in S$ we have $s \cdot b \in S$.*

On similar lines we can define Smarandache definite left ideals.

We illustrate this by the following example.

***Example 1.2.28:*** Let $Z_{12} = \{\bar{0}, \bar{1}, \bar{2}, \ldots, \overline{11}\}$ be the ring of integers modulo 12. Let $A = \{\bar{0}, \bar{4}, \bar{8}\}$ be a field. $Z_{12}$ is a S-ring $S = \{0, 6\}$ is a S-definite ideal related to A. But S is also an ideal $Z_{12}$.

The interested reader is requested to refer [35-8].
  Now we proceed on to recall the definition of semirings.

**DEFINITION 1.2.10:** *Let S be a non empty set on which is defined two binary operations addition '+' and multiplication '∘' satisfying the following conditions.*
1. *$(S, +)$ is a commutative monoid.*
2. *$(S, \circ)$ is a semigroup.*
3. *$(a + b) \circ c = a \circ c + b \circ c$ and*

*$c \circ (a + b) = c \circ a + c \circ b$ for all a, b, c in S.*

  *We call $(S, +, \circ)$ to be a semiring. If the semiring $(S, +, \circ)$ contains 1 such that $1 \circ s = s \circ 1 = s$ for all $s \in S$, we call S to be a semiring with unit.*

We illustrate this by the following examples.

**Example 1.2.29:** Let $Z^o = Z^+ \cup \{0\}$ the set of integers with zero $(Z^o, +, \circ)$ is a semiring.

**Example 1.2.30:** Let $Q^o = Q^+ \cup \{0\}$; $Q^o$ is a commutative semiring.

*Note:* A semiring S is said to be commutative if $a \circ b = b \circ a$ for all $a, b \in S$.

**Example 1.2.31:** Let
$$S_{2\times 2} = \left\{ \begin{pmatrix} a & b \\ c & d \end{pmatrix} \middle| a, b, c, d \in Z^o = Z^+ \cup \{0\} \right\}$$
set of all $2 \times 2$ matrices with entries from $Z^o$. $(M_{2\times 2}, +, \circ)$ is a semiring '+' is the matrix addition and '$\circ$' is the matrix multiplication.

We now proceed on to define the notion of subsemiring.

**DEFINITION 1.2.11:** *Let S be a semiring, P a proper subset of S. We call P to be subsemiring if P itself is a semiring, under the operations of S.*

We illustrate this situation by the following examples.

**Example 1.2.32:** Let $Z^o = Z^+ \cup \{0\}$ be a semiring. Take $2Z^o = \{0, 2, 4, \ldots\}$; $2Z^o$ is a subsemiring of $Z^o$.

**Example 1.2.33:** Let $Z^o[x]$ be the polynomial semiring i.e., $Z^o[x] = \{p_0 + p_1 x + \ldots + p_n x^n \mid p_0, p_1, \ldots, p_n \in Z^o\}$. Under the polynomial addition and multiplication $Z^o[x]$ is a semiring. Clearly $Z^o \subset Z^o[x]$ and $Z^o$ is the subsemiring of $Z^o[x]$.

Now as in case of rings in case of semirings also one can define the notion of ideals of a semiring. We now recall the definition of ideals of a semiring and illustrate it with examples. For more about semirings the interested reader is requested to refer [92, 96, 98, 102].

**DEFINITION 1.2.12:** *Let S be a semiring. A non empty subset I of S is said to be an ideal of S if*
 1. *I is a subsemiring of S.*
 2. *For all $i \in I$ and $s \in S$ we have $i\,s$ and $s\,i \in I$.*

As in case of rings one can define both right ideal of a semiring and a left ideal of a semiring. The notion of right ideal and left ideal trivially coincide only when the semiring is a commutative semiring.

***Example 1.2.34:*** Let $Z^o = Z^+ \cup \{0\}$ be a semiring n $Z^o$, n any positive integer is an ideal of the semiring $Z^o$. We have n $Z^o$ = $\{0, n, 2n, 3n, \ldots\}$.

***Example 1.2.35:*** Consider the distributive lattice given by the following figure.

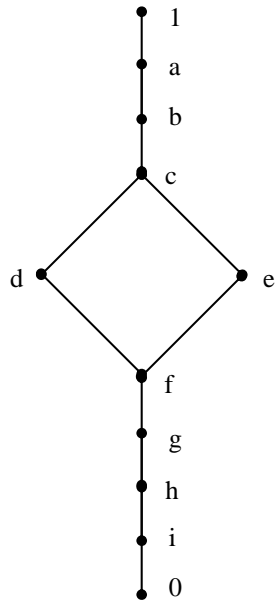

**Figure 1.2.1**

S = {1, a, b, c, d, e, f, g, h, i, 0}. S is a semiring. Take I = {0, i, h, g, f, d, e, c} is an ideal of the semiring S.

As in case of ring in case of semirings one can define the notion of units, idempotents and zero divisors; but in case of semirings the existence of idempotents need not lead to zero divisors.

For we see all distributive lattices are semirings and every element in a distributive lattice is an idempotent but this does not imply the semiring has zero divisors.

We illustrate this by the following example.

*Example 1.2.36:* Consider semiring given by the lattice.

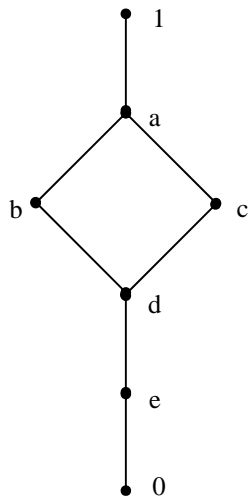

**Figure 1.2.2**

S = {1, a, b, c, d, e, 0} is a semiring but S has no zero divisors but every element in S is an idempotent of S.

But we have semiring which contain zero divisors as well as idempotents.

*Example 1.2.37:* Consider the semiring given by the following lattice.

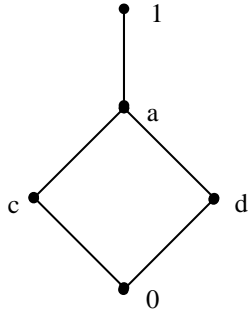

**Figure 1.2.3**

$S = \{1, a, c, d, 0\}$; S is a semiring in which every element; is an idempotent. This semiring has zero divisors. For $c \cdot d = 0$; $c \neq 0$ and $d \neq 0$. Hence the claim.

We give yet another example.

*Example 1.2.38:* Consider the semiring

$$M_{2\times 2} = \left\{ \begin{pmatrix} a & b \\ c & d \end{pmatrix} \middle| a, b, c, d \in Z^o = Z^+ \cup \{0\} \right\}.$$

$M_{2\times 2}$ is a semiring under matrix addition and matrix multiplication.
    Take

$$A = \begin{pmatrix} 5 & 0 \\ 0 & 0 \end{pmatrix} \text{ and } B = \begin{pmatrix} 0 & 0 \\ 0 & 7 \end{pmatrix}.$$

Clearly

$$A \circ B = \begin{pmatrix} 5 & 0 \\ 0 & 0 \end{pmatrix} \begin{pmatrix} 0 & 0 \\ 0 & 7 \end{pmatrix} = \begin{pmatrix} 0 & 0 \\ 0 & 0 \end{pmatrix}.$$

Thus the semiring $M_{2\times 2}$ has zero divisors.

Now we proceed on to define the notion of semifields.

**DEFINITION 1.2.13:** *Let S be a non empty set. S is said to be a semifield if*

1. *S is a commutative semiring with 1.*
2. *S is a strict semiring i.e., for a, b ∈ S if a + b = 0 then a = 0 and b = 0.*
3. *If in S, a ∘ b = 0 then either a = 0 or b = 0.*

We now illustrate this definition by the following examples.

***Example 1.2.39:*** $Z^o = Z^+ \cup \{0\}$ be a semiring which is also a semifield.

***Example 1.2.40:*** The semiring

$$M_{2\times 2} = \left\{ \begin{pmatrix} a & b \\ c & d \end{pmatrix} \middle| a, b, c, d \in Z^o = Z^+ \cup \{0\} \right\}$$

is not a semifield for $M_{2\times 2}$ is not a commutative semiring and further we have

$$A = \begin{pmatrix} 3 & 0 \\ 0 & 0 \end{pmatrix} \text{ and } B = \begin{pmatrix} 0 & 0 \\ 0 & 2 \end{pmatrix}$$

such that

$$A \circ B = \begin{pmatrix} 3 & 0 \\ 0 & 0 \end{pmatrix} \begin{pmatrix} 0 & 0 \\ 0 & 2 \end{pmatrix} = \begin{pmatrix} 0 & 0 \\ 0 & 0 \end{pmatrix}.$$

i.e., $M_{2\times 2}$ has zero divisors as A o B = (0) without A being zero and B being equal to zero.

Now we proceed on to define the notion of subsemifield.

**DEFINITION 1.2.14:** *Let S be a semifield; a proper subset P of S is said to be subsemifield if P itself is a semifield.*

***Example 1.2.41:*** Let $Q^o = Q^+ \cup \{0\}$ be the semifield. $Z^o \subset Q^o$ and $Z^o$ is a subsemifield of $Q^o$.

We see all semifields need not in general contain sub semifields. For instance $Z^o$ has no proper subset which is a sub semifield of $Q^o$.

**Example 1.2.42:** We see $Z^o[x]$ is a semifield. $Z^o \subset Z^o[x]$ is a subsemifield.

Now in the next section we proceed on to define the notion of vector spaces, semivector spaces and modules.

## 1.3 Vector spaces, Semivector spaces and Modules

In this section we just recall the basic definitions of vector spaces, modules and semivector spaces and illustrate them with examples. For more about these notions please refer [2, 15, 29, 104].

**DEFINITION 1.3.1:** *We call a non empty set V to be a vector space over a field F of characteristic zero if the following conditions are satisfied.*

1. *$(V, +)$ is an additive abelian group.*
2. *For every $v \in V$ and $a \in F$ (The elements of the field F are called as scalars) we have $v \cdot a$ and $a \cdot v \in F$. (Note $a \cdot v = v \cdot a$ we always write the scalar as coefficients).*
3. *$0 \cdot v = v \cdot 0 = 0$ for all $v \in V$ and $0 \in F$.*
4. *$1 \cdot v = v \cdot 1 = v$ for all $v \in V$.*
5. *$(a + b) v = av + bv$ for all $a, b \in F$ and $v \in V$.*
6. *$a (u_1 + u_2) = a u_1 + a u_2$ for all $a \in F$ and $u_1, u_2 \in V$.*

*Since $a (u_1 + u_2) = (u_1 + u_2) a$ and $(a + b) v = v (a + b)$, we do not express both identities explicitly.*

We illustrate this definition by an example.

**Example 1.3.1:** Let $Q[x]$ be the polynomial abelian group with respect to addition. Q the field of rationals. $Q[x]$ is a vector space over Q.

*Note:* Q[x] is not a vector space over the field of reals. For if we take $\sqrt{2} \in R$ and $p(x) = p_0 + p_1 x + \ldots + p_n x^n$ any polynomial in Q[x] where $p_0, p_1, \ldots, p_n \in Q$. We see $\sqrt{2}\, p(x)$ is not a polynomial in Q[x], hence the claim.

Now we show more interesting result in this direction.

*Example 1.3.2:* Let R[x] be the set of all polynomials with coefficients from the reals in the indeterminate x. R[x] is a group under addition. R[x] is a vector space over Q as well as R[x] is a vector space over R.

*Example 1.3.3:* Let R be the group of reals under addition. R is a vector space over R as well as R is a vector space over Q.

*Example 1.3.4:* Let Q be the group of rationals under addition. Q is a vector space over the field of rationals Q; but Q is not a vector space over the field of reals. For $3 \in Q$ and $\sqrt{5} \in R$ but $3\sqrt{5} \notin Q$, hence the claim.

*Example 1.3.5:* Let $M_{3\times 2}$ be the collection of all $3 \times 2$ matrices with entries from the field of rationals i.e., $M_{3\times 2} = M = (m_{ij})$ is a group under matrix addition. $M_{3\times 2}$ is a vector space over the field Q but $M_{3\times 2}$ is not a vector space over the field of reals R for take $\sqrt{3} \in R$, clearly

$$\sqrt{3}M = \begin{pmatrix} \sqrt{3}\, m_{11} & \sqrt{3}\, m_{12} \\ \sqrt{3}\, m_{21} & \sqrt{3}\, m_{22} \\ \sqrt{3}\, m_{31} & \sqrt{3}\, m_{32} \end{pmatrix};$$

but $\sqrt{3}$ M does not belong to the collection $M_{3\times 2}$ as $\sqrt{3} m_{ij} \notin Q$. Hence the claim. Clearly $M_{3\times 2}$ is a vector space over Q only.

Now we have seen examples of vector spaces. We recall the definition of linear algebra.

**DEFINITION 1.3.2:** *Let V be a vector space over the field F, i.e., V is an additive abelian group satisfying the conditions given in definition 1.3.1. If in addition we have for u, v ∈ V, v · u and u · v ∈ V and the operation '·' is a closed binary associative operation on V. Then we call V to be linear algebra over the field F.*

We see all linear algebras are vector spaces over the fields on which they are defined, but in general a vector space need not always be a linear algebra. For more refer [2, 15, 29, 104].

We now illustrate this definition.

*Example 1.3.6:* Let $M_{3\times 2}$ be a vector space over Q given in example 1.3.5. Clearly $M_{3\times 2}$ is not a linear algebra over Q as for any A, B ∈ $M_{3\times 2}$ we see the matrix product AB is not even defined.

*Example 1.3.7:* Let Q[x] be the group of polynomials under addition. Q[x] is a vector space over Q. In fact Q[x] is a linear algebra over Q; for if p(x), q(x) ∈ Q[x] then p(x)q(x) the product of two polynomials is again polynomial in Q[x]. Hence Q[x] is a linear algebra over Q.

*Example 1.3.8:* Let $M_{n\times n}$ = {M = $(m_{ij})|m_{ij}$ ∈ Q} be the collection of n × n matrices with entries from Q. $M_{n\times n}$ is a vector space over Q. In fact $M_{n\times n}$ is a linear algebra over Q. In fact $M_{n\times n}$ is neither a vector space nor a linear algebra over R.

We have seen the notion of vector spaces and linear algebra. Now we just define the basis for a vector space and dimension of a vector space over a field K.

**DEFINITION 1.3.3:** *Let V be a vector space over the field K. Let B = { $v_1$, $v_2$, …, $v_n$} be a non empty subset of V; we say $v_1$, $v_2$, …,* $v_n$ are linearly dependent over K if there exits elements $\alpha_1, \alpha_2, ..., \alpha_n$ in K not all of them zero, such that $\alpha_1 v_1 + \alpha_2 v_2 + ... + \alpha_n v_n = 0$. If the vectors $v_1, v_2, ..., v_n$ are not linearly dependent over K then they are said to be linearly independent; that is

$$\sum_{i=1}^{n} \alpha_i v_i = 0$$

if and only if each $\alpha_i = 0$ for $i = 1, 2, ..., n$.

We say the set $B = \{v_1, v_2, ..., v_n\}$ which is linearly independent generates the vector space V over K if every element x in V can be uniquely represented as a linear combination of elements from B with coefficients from K i.e.,

$$x = \sum_{i=1}^{n} x_i v_i$$

where $x_i \in K; 1, 2, ..., n$.

The set B in this case is called a basis of V and the number of elements in B corresponds to the dimension of V over K. If B has finite number of elements in it then we call V to be a finite dimensional vector space over K. If B has infinite number of elements then we say V is an infinite dimensional vector space over K.

It is important to note that the dimension is also dependent on the field over which it is defined.

This will be clear from the following examples.

*Example 1.3.9:* Let $M_{2 \times 2} = \{(m_{ij}) \mid m_{ij} \in R\}$ be set of all $2 \times 2$ matrices with entries from the field of reals. $M_{2 \times 2}$ is a vector space over R. Now the dimension of $M_{2 \times 2}$ as a vector space over R is four i.e., $M_{2 \times 2}$ is generated by the set

$$B = \left\{ \begin{pmatrix} 1 & 0 \\ 0 & 0 \end{pmatrix}, \begin{pmatrix} 0 & 1 \\ 0 & 0 \end{pmatrix}, \begin{pmatrix} 0 & 0 \\ 1 & 0 \end{pmatrix}, \begin{pmatrix} 0 & 0 \\ 0 & 1 \end{pmatrix} \right\}.$$

Now $M_{2\times 2}$ is also a vector space over Q; but the dimension of $M_{2\times 2}$ over Q is infinite dimensional.

This can be more explicitly explained by the following example.

***Example 1.3.10:*** Let $V = R \times R$ be a vector space over R. V is a two dimensional vector space over the reals. V is generated by $B = \{(0, 1), (1, 0)\}$. Now $V = R \times R$ is also a vector space over Q.

Now V as a vector space is generated by an infinite number of elements as the generating pairs $B = \{(x, y)\}$ must vary over the irrationals of the form $(\sqrt{p}, \sqrt{q})$ where p and q are not perfect squares i.e.,

$$B = \{(\sqrt{2}, 0), (0, \sqrt{2}), (\sqrt{3}, 0), (0, \sqrt{3}), (\sqrt{5}, 0), (0, \sqrt{5}), ...\}.$$

Thus these two examples clearly show the dimension of a vector space is much dependent on the field over which it is defined.

Now we proceed on to define the notion of modules.

**DEFINITION 1.3.4:** *Let R be a ring. A non empty set M is said to be a R-module (or a module over R) if M is an abelian group under an operation '+' such that for every $r \in R$ and $m \in M$ there exists an element rm in M subject to*

$$r(a + b) = ra + rb$$
$$r(sa) = (rs)a$$
$$(r + s)a = ra + sa$$

*for all $a, b \in M$ and $r, s \in R$.*

*If the R has a unit element 1 and if $1 \cdot m = m$ for every element m in M; then M is called a unital R-module.*

Note that if R is a field, a unital R-module is nothing more than a vector space over R.

We will be defining and using only unital rings. Properly speaking, we should have called the object that we have defined as a left R-module. An R-module M is said to be cyclic if there is an element $m_o \in M$ such that every $m \in M$ is of the form $m = rm_o$ where $r \in R$.

An additive subgroup A of the R-module M is called a submodule of M if whenever $r \in R$ and $a \in A$ then $ra \in A$.

An R-module M is said to be finitely generated if there exists elements $a_1, a_2, \ldots, a_n \in M$ such that every m in M is of the form $m = r_1a_1 + \ldots + r_na_n$.

We now illustrate this by the following example.

*Example 1.3.11:* Let Z be the ring of integers

$$M_{2\times 2} = \left\{ \begin{pmatrix} a & b \\ c & d \end{pmatrix} \middle| a, b, c, d \in Q \right\}$$

is an additive group. $M_{2\times 2}$ is a module over Z.

*Example 1.3.12:* Let Z[x] be the polynomial ring. Z[x] is a unital module over Z.

*Example 1.3.13:* Let $M_{5\times 3} = \{(m_{ij}) \mid m_i \in Q\}$ be the set of all $5 \times 3$ matrices with entries from Q. $M_{5\times 3}$ is a group under matrix addition. $M_{5\times 3}$ is a vector space over Q also $M_{5\times 3}$ is a vector space over Z.

Let Q be a field $Q(\sqrt{3})$ is a field containing in Q. $Q(\sqrt{3})$ is an extension field of Q.

$Q(\sqrt{3})$ is a vector space over Q and in fact $Q(\sqrt{3})$ is a finite dimension vector space over Q; then we call $Q(\sqrt{3})$ a finite extension of Q and $\{1, \sqrt{3}\}$ is a basis and $Q(\sqrt{3})$ is a vector space of dimension 2.

Thus in view of this we have the following definition.

**DEFINITION 1.3.5:** *Let F be any field, E an extension field of F that is F is a subfield of the field E. We say E is a finite extension of F if E as a vector space over F is finite dimensional and it is denoted by [E: F] = n < ∞.*

We see F can be realised as a subfield of E. One of the interesting properties is that if F is any field; E is a finite extensions of F and T is a finite extension of E then T is a finite extension of F i.e., F ⊆ E ⊆ T; i.e., [T: F] = [T: E] [E: F] = mn where [T: E] = m < ∞ and [E: F] = n < ∞.

If T is a vector space over E of dimension m and E is a vector space of dimension n then T is a vector space of dimension mn. Several interesting results can be obtained in this direction but the reader is requested to refer [2, 15, 29, 104] for more about this literature.

We see Q is a field and R the field of reals is only an infinite extension as R is a vector space over Q but dimension of R as a vector space over Q is infinite. Also if Q is the field of rationals $Q(\sqrt{2})$ is a finite extension of Q.

Now consider
$$K = Q(\sqrt{2}),$$
now $K(\sqrt{3})$ is a finite extension of K. We see $Q(\sqrt{2})$ is a vector space of dimensions two over Q and $K(\sqrt{3})$ is a vector space of dimension two over K. We see
$$Q \subseteq K \subseteq K(\sqrt{3}).$$

Now $K(\sqrt{3})$ is a vector space over Q, and $K(\sqrt{3})$ is a finite dimensional vector space over Q.

Now

$$\left[ K\left(\sqrt{3}\right) : Q \right] = \left[ K\left(\sqrt{3}\right) : K \right] [K : Q] = 2.2 = 4.$$

Thus $K\left(\sqrt{3}\right)$ is a vector space of dimension 4 over Q and a basis of $K\left(\sqrt{3}\right)$ as a vector space over Q is $\{1, \sqrt{2}, \sqrt{3}, \sqrt{6}\}$.

Now having seen finite extensions of fields, we now proceed on to recall the definition of Smarandache vector space.

For more about Smarandache vector space refer [104].

**DEFINITION 1.3.6:** *Let R be a S-ring, V be a module over R. We say V is a Smarandache vector space of type II (S-vector space of type II) if V is a vector space over a proper subset K of R where K is a field.*

We illustrate this by the following example.

*Example 1.3.14:* R[x] be the polynomial ring over the field of reals. Q is a field in R[x]. Q[x] is a S-vector space over Q. Clearly Q[x] is not a vector space over R.

**DEFINITION 1.3.7:** *Let R be a S-ring M a R-module M is said to be a Smarandache linear algebra of type II (S-linear algebra of type II) if M is a linear algebra over a proper subset K in R where K is a field.*

For more about these concepts please refer [101].

**DEFINITION 1.3.8:** *Let M be an R-module over a S-ring R. If a proper subset P of M is such that P is a S-vector space II over a proper subset K of R where K is a field then we call P a Smarandache subspace II(S-subspace II) of M relative to P.*

Now we proceed on to define the notion of Smarandache linear subalgebra.

**DEFINITION 1.3.9:** *Let M be an R-module over a S-ring R. If M is a proper subset P such that P is a Smarandache linear algebra over a proper subset K in R where K is a field then we say P is a Smarandache linear subalgebra of type II (S-linear subalgebra of type II).*

For more refer [104].

**Chapter Two**

# SMARANDACHE ALGEBRAIC SPECIAL DEFINITE STRUCTURES

The study of Smarandache algebraic structures was only visualizing an algebraic stronger structure in an algebraically weak structure. For instance the Smarandache semigroups had in it a group, a Smarandache groupoid had a loop in it as a substructure and so on. It is pertinent to mention here that we had no concept of Smarandache group for the group happened to be a very strong concentrated algebraic structure for no more stronger structure can be put into that algebraic structure.

So now for the first time we study the new notion of Smarandache special definite structures. That is we visualize a weak algebraic structure in a strong algebraic structure.

This chapter has two sections. In section one the notion of Smarandache special definite groups are introduced. Smarandache special definite algebraic structures or equivalently we say Smarandache special definite algebraic structure like S-special definite fields, S-definite special minimal (and maximal ideals) S-special definite homomorphism are introduced.

## 2.1 Smarandache Special Definite Groups

In this section we for the first time introduce the notion of Smarandache special definite groups and investigate some of its interesting properties.

**DEFINITION 2.1.1:** *Let (G, \*) be a group we call G to be Smarandache special definite group (S-special definite group) if we can find a non empty subset S in G such that (S, \*) is a semigroup.*

We illustrate this situation by the following examples.

*Example 2.1.1:* Let $(Z, +)$ be a group under addition. If we take $Z^+$ the set of only positive integers then $Z^+$ is a semigroup under addition. Thus $(Z, +)$ is a Smarandache special definite group.

*Example 2.1.2:* Consider $G = Q \setminus \{0\}$; G is a group under multiplication ×; Take $Z \setminus \{0\}$ a proper subset of G, clearly $Z \setminus \{0\}$ is only a semigroup under multiplication. Thus $(Q \setminus \{0\}, \times)$ is a S-special definite group.

All groups in general are not Smarandache special definite groups. This is proved only by the following example.

*Example 2.1.3:* Consider $S_3$ the group of permutations on (1 2 3). Clearly $S_3$ is not a S-special definite group for we cannot find any proper subset in $S_3$ which is a semigroup under the operations of $S_3$.

*Note:* When we say a proper subset $(S, *) \subseteq (G, *)$ is a semigroup we by no means accept $(S, *)$ to be a subgroup of $(G, *)$.

*Example 2.1.4:* $S_n$ ($n \geq 2$) is a symmetric group which is never a S-special definite group.

***Example 2.1.5:*** Let $M_{2\times 2} = \{(a_{ij}) = M$ such that $a_{ij} \in Q$ and $|M| \neq 0$, i.e., determinant of $M \neq 0\}$ be a group under matrix multiplication. Take $P_{2\times 2} = \{(a_{ij}) = P \mid a_{ij} \in Z, |M| \neq 0\}$ a proper subset of $M_{2\times 2}$. Clearly P is only a semigroup under multiplication. So $M_{2 \times 2}$ is a S-special definite group.

We make a note $M_{2\times 2}$ is a non commutative group.

***Example 2.1.6:*** Consider the set $M_{2\times 2} = \{(m_{ij}) = M \mid m_{ij} \in Q$ and $\det[M] \neq (0)\}$. $M_{2\times 2}$ is a group under matrix multiplication. Take

$$P_{2\times 2} = \left\langle \begin{pmatrix} 1 & -2 \\ -2 & 1 \end{pmatrix} \right\rangle ;$$

$P_{2\times 2}$ is a semigroup under matrix multiplication and $P_{2\times 2}$ is a proper subset of $M_{2\times 2}$. So $M_{2 \times 2}$ is a S-special definite group.

**DEFINITION 2.1.2:** *Let (G, \*) be a group. If we have a proper subset P in G such that (P, \*) is a commutative semigroup then we call (G, \*) to be a commutative Smarandache special definite group (commutative S-special definite group).*

*Note:* It is important and interesting to note that (G, \*) need not be commutative group still the group (G, \*) may be a commutative S-special definite group.

We have got the following theorem as an immediate consequence of the definition which is left as an exercise for the reader.

**THEOREM 2.1.1:** *Every commutative group (G, \*) which is a S-definite special group is a commutative S-definite special group.*

Just we show by an example that there exists commutative S-special definite groups (G, \*) which in general are not commutative.

***Example 2.1.7:*** Let $M_{2\times 2} = \{(m_{ij}) = M \mid m_{ij} \in Q$ and $|M| \neq 0\}$ be the group under matrix multiplication.

$$P_{2\times 2} = \left\langle \begin{pmatrix} 1 & -2 \\ -2 & 1 \end{pmatrix} \right\rangle$$

is a semigroup generated under matrix multiplication. We see $P_{2\times 2}$ is a commutative semigroup but $M_{2\times 2}$ is not a commutative group but $M_{2\times 2}$ is a commutative S-special definite group.

Now we proceed on to define strongly commutative Smarandache special definite group.

**DEFINITION 2.1.3:** *If every proper subset P of G which is a semigroup under "*" happens to be a commutative semigroup then we call (G, *) to be strongly commutative Smarandache special definite group (strongly commutative S-special definite group).*

The following theorem is left as an exercise for the interested reader to prove.

**THEOREM 2.1.2:** *Let (G, *) be a commutative group, which is also a S-special definite group. Then (G, *) is a strongly commutative Smarandache special definite group.*

Now we proceed on to prove that all commutative S-special definite groups are not strongly commutative S-special definite groups but all strongly commutative S-special definite groups are commutative S-special definite groups. It is clear from the very definition all strongly commutative S-special definite groups are commutative S-special definite groups.

We prove by an example that a commutative S-special definite group in general is not a strongly commutative S-special definite group.

*Example 2.1.8:* Consider $M_{2\times 2} = \{(M = m_{ij}) \mid m_{ij} \in Q$ and $|M| \neq 0\}$ be a non commutative group under matrix multiplication. Consider the set $P_{2\times 2} = \{P = (p_{ij}) \mid p_{ij} \in Z, |P| \neq 0\}$, $P_{2\times 2}$ is a subset of $M_{2\times 2}$ and $P_{2\times 2}$ is a semigroup under matrix multiplication and $P_{2\times 2}$ is not a commutative semigroup.
Now consider

$$T_{2\times 2} = \left\langle \begin{pmatrix} 1 & -2 \\ -2 & 1 \end{pmatrix} \right\rangle$$

i.e., $T_{2 \times 2}$ is a semigroup generated by the matrix

$$\begin{pmatrix} 1 & -2 \\ -2 & 1 \end{pmatrix}$$

which is in $M_{2\times 2}$. Clearly $T_{2\times 2}$ is commutative and so $M_{2\times 2}$ is only a commutative S-special definite group and not a strongly commutative S-special definite group of a S-special definite group.

**DEFINITION 2.1.4:** *Let (G, \*) be a group. Let (H, \*) be a subgroup of (G, \*) If (H, \*) is itself a S-special definite group then we call (H, \*) to be a S-special definite subgroup of (G, \*).*

We have the following theorem.

**THEOREM 2.1.3:** *Let (G, \*) be a group. Suppose (H, \*) be a subgroup of (G, \*) which is a S-special definite subgroup. Then (G, \*) is a S-special definite group.*

*Proof:* Given (G, \*) is a group and H is a subgroup of G which is a S-special definite group. So (H, \*) has a proper subset P such that (P, \*) is a semigroup of H. Now $P \subset H$ and $H \subset G$ so $P \subset H \subset G$ and H is a semigroup, hence G is a S-special definite group.
  Hence the claim.

We illustrate this situation by the following example.

*Example 2.1.9:* Let $Q \setminus \{0\}$ be the group under multiplication. Consider $Q^+$ the set of positive rationals, $Q^+$ is a subgroup of $Q \setminus \{0\}$ under multiplication. Let $Z^+$ be the semigroup under multiplication contained in $Q^+$. So $Q^+$ is a S-definite special subgroup of Q. Since $Z^+ \subset Q^+ \subset Q \setminus \{0\}$ we see $Q \setminus \{0\}$ is also a S-definite special group.
  We see all subgroups of a S-definite special group need not be a S-definite special subgroup.

We illustrate this by the following example.

***Example 2.1.10:*** Let $G = S_3 \times (Q \setminus \{0\}, \times)$ be a group. Clearly G is a S-definite special group for $P = \{e\} \times (Z^+, \times)$ is a semigroup of G, so the claim.

Now we see $A_3 \times \{1\}$ is also a subgroup of G but $A_3 \times \{1\}$ is not a S-definite special subgroup of G. $M = \{e\} \times (Q^+, \times)$ is a S-definite special subgroup of G. Hence the claim.

Now we proceed on to define the notion of S-special definite normal subgroup of a group G.

**DEFINITION 2.1.5:** *Let (G, \*) be a group (H, \*) be a normal subgroup of (G, \*), we call (H, \*) to be a Smarandache definite special normal subgroup (S-definite special normal subgroup) of G if (H, \*) is itself a S-special definite group. If (G, \*) has no S-special definite normal subgroups but (G, \*) is a S-special definite group then we call (G, \*) to be a Smarandache definite special simple group (S-definite special simple group).*

It is interesting to make a note of the fact that if (G, \*) is a commutative S-definite special group having non trivial S-definite special subgroup then G has S-definite special normal subgroups, evident from the fact every subgroup of a commutative group is normal.

**THEOREM 2.1.4:** *Let (G, \*) be a group having a S-definite special normal subgroup. Then (G, \*) is a S-special definite group.*

The reader is expected to prove this theorem which is similar to the proof of the theorem 2.1.3.

It is important to mention here that torsion free groups alone happen to be S-special definite groups or at least the group should contain some torsion free elements otherwise they by Cauchy theorem turn to be not a S-special definite groups. This condition has sealed the opportunity of us to study the classical theorems like Sylow, Lagrange and Cauchy.

Now we define element wise Smarandache properties of these S-special definite groups.

**DEFINITION 2.1.6:** *Let (G, \*) be a S-special definite group. A proper subset P of (G, \*) is said to be Smarandache special definite ideal (S-special definite ideal) of (G, \*) if P is an ideal with respect to some semigroup T of G; i.e., tp, pt $\in$ P for all t $\in$ T and p $\in$ P.*

It is interesting and important to note that every semigroup of the group (G, \*) need not have an ideal associated with it or equivalently the S-special definite ideal of G need not be an ideal with every semigroup of (G, \*) if it is an ideal with respect to some semigroup it is enough.

We illustrate this by the following example.

*Example 2.1.11:* Let Q \ {0} be the group under multiplication. $Z^+$ is a subset of Q \ {0} which is also a semigroup. Thus (Q \ {0}, ×) is a S-special definite group. Now $2Z^+$ is an ideal of the S-special definite semigroup of Q \ {0}. Hence $2Z^+$ is S-special definite ideal of $Z^+$. Take T = Z \ {0} $\subseteq$ Q \ {0}; T is also a S-special definite ideal of Q \ {0} for T is a ideal over $Z^+$.
    Further we see $Z^+$ is not a S-special definite ideal of Q \ {0}.

We define the notion of Smarandache special definite cyclic ideal, Smarandache special definite maximal ideal, Smarandache special definite minimal ideal and Smarandache special definite prime ideal related with a S-special definite group.

**DEFINITION 2.1.7:** *Let (G, \*) be a S-special definite group i.e., (G, \*) has at least a semigroup (T, \*) such that T $\subseteq$ G. We say (P, \*) is a Smarandache special definite maximal ideal (S-special definite maximal ideal) of (G, \*) related to T if (M, \*) is any other S-special definite ideal related to T and if P $\subseteq$ M $\subseteq$ T then either P = M or M = T. We say (U, \*) to be a S-minimal definite ideal related to T if we have another (V, \*) related to T and (0) = V $\subseteq$ U $\subseteq$ T then (0) = V or V = U. We call an S-* *special definite ideal W related to T to be a Smarandache special definite prime ideal (S-special definite prime ideal) if for ab ∈ W either a or b is in W.*

Now having defined these new Smarandache notions we proceed on to describe them with examples.

*Example 2.1.12:* Let $Q \setminus \{0\}$ be the group under multiplication. $Q \setminus \{0\}$ is a S-special definite group. Let $Z \setminus \{0\} \subseteq Q \setminus \{0\}$ be the semigroup of $Q \setminus \{0\}$ under multiplication. $2Z \setminus \{0\}$ is a S-definite special ideal of Z. Clearly $2Z \setminus \{0\}$ is S-definite special maximal ideal of $Z \setminus \{0\}$. All S-definite ideals of $Z \setminus \{0\}$ are not maximal for example; consider $6Z \setminus \{0\}$ a proper subset of $Z \setminus \{0\}$. $6Z \setminus \{0\}$ is a S-special definite ideal related to the semigroup $Z \setminus \{0\}$, but $6Z \setminus \{0\}$ is not S-definite special maximal ideal of $Z \setminus \{0\}$ for the S-definite special ideal $2Z \setminus \{0\}$ contains $6Z \setminus \{0\}$ to be a proper subset. Hence the claim. Take $2Z \setminus \{0\}$ the S-definite special ideal of $Z \setminus \{0\}$; $2Z \setminus \{0\}$ is a S-special definite prime ideal of $Z \setminus \{0\}$.

We are not able to get any S-special definite minimal ideals of a group, we leave it as an open problem in the last chapter. Now having seen illustrations of some types of S-special definite ideals we define now the notion of Smarandache special definite principal or cyclic ideals related to the semigroup of a S-special definite group.

**DEFINITION 2.1.8:** *Let (G, \*) be a S-special definite group. Suppose T be a proper subset of G which is a semigroup. We call P a proper subset of T to be a Smarandache special definite principal ideal (S-special definite principal ideal) related to P if P is generated by a single element.*

*Example 2.1.13:* Let $Q \setminus \{0\}$ be a S-special definite group. $Z \setminus \{0\}$ be a semigroup of $Q \setminus \{0\}$. $6Z \setminus \{0\}$ is an S-special definite ideal of $Z \setminus \{0\}$, we see $6Z \setminus \{0\}$ is a S-special definite principal ideal of $Z \setminus \{0\}$.

*Remark:* It is interesting to note that in general a S-special definite principal ideal need not be a S-special definite maximal ideal related to some semigroup of a S-special definite group G.

We suggest some problems in the last chapter of this book. All S-special definite principal ideals related to a semigroup of a S-special definite group need not be prime.
    We illustrate this by the following example.

*Example 2.1.14:* Let $Q \setminus \{0\}$ be a S-special definite group under multiplication. $Z \setminus \{0\}$ be a semigroup of $Q \setminus \{0\}$. Consider $6Z \setminus \{0\}$ the S-special definite principal ideal related to $Z \setminus \{0\}$ of the S-special definite group $Q \setminus \{0\}$. $6Z \setminus \{0\}$ is generated by 6. Now $4.3 \in 6Z \setminus \{0\}$ but neither 4 nor 3 is in $6Z \setminus \{0\}$. Thus a S-special definite principal ideal $6Z \setminus \{0\}$ is not a S-special definite prime ideal of $Z \setminus \{0\}$. Hence the claim.

Now having defined the S-special definite ideals we proceed on to define Smarandache definite coset of a S-special definite ideals of a S-special definite group.

**DEFINITION 2.1.9:** *Let (G, \*) be a S-special definite group. H be a semigroup of G. We define the Smarandache definite left coset of the semigroup (S-definite left coset of the semigroup ) aH of H in (G, \*) for a $\in$ G as follows.*
$$aH = \{ah \mid h \in H \}.$$
    *If a $\in$ H then aH $\neq$ H in general. Similarly Smarandache definite right coset of the semigroup (S-definite right coset of the semigroup) H is defined to be Ha of H in (G, \*) for a $\in$ G as follows:*
$$Ha = \{ha \mid h \in H\}.$$

*Remark:* If (G, \*) be a S-special definite group and if (G, \*) is a commutative group then we have the S-definite special right coset to be equal to S-definite special left coset i.e., Ha = aH.

We illustrate this situation by the following example.

***Example 2.1.15:*** Let $Q \setminus \{0\}$ be a S-special definite group. $Z^+ \subset Q \setminus \{0\}$ be a semigroup $g Z^+ = Z^+ g$ for all $g \in Q \setminus \{0\}$, $gZ^+$ where $g = \dfrac{1}{2}$ then

$$\frac{1}{2}Z^+ = \left\{\frac{1}{2}, 1, \frac{3}{2}, 2, \frac{5}{2}, 3,\ldots\right\}.$$

Also for $2 \in Z^+$; $2Z^+ \neq Z^+$.

**THEOREM 2.1.5:** *Let G be a S-special definite group. H be a semigroup of G. Let aH be the S-definite coset of the semigroup H in G. Then aH $\neq$ H for a $\in$ H. We may have H $\subseteq$ gH for g $\in$ G \ H.*

*Proof:* We prove this only by a counter example. Let $Q \setminus \{0\}$ be the S-special definite group under multiplication. Take $Z^+ \subset Q \setminus \{0\}$ to be the semigroup of $Q \setminus \{0\}$.

Take $\dfrac{1}{2} \in Q \setminus \{0\}$ now

$$\frac{1}{2}Z^+ = \left\{\frac{1}{2}, 1, \frac{3}{2}, 2, \frac{5}{2}, 3, \frac{7}{2}, 4, \frac{9}{2}, 5, \frac{11}{2}, 6, \ldots, \infty\right\}.$$

Clearly $Z^+ \subset \dfrac{1}{2}Z^+$. Hence aH may contain H. Consider $2 \in Z^+$, the coset $2Z^+ = \{2, 4, 6, \ldots\} \subset Z^+$. Thus even if $a \in H$, $aH \neq H$. Hence the claim.

This is the marked difference between the coset of a subgroup and S-definite coset of a semigroup.

**THEOREM 2.1.6:** *Let (G, \*) be a S-special definite group. H a semigroup of G. The S-special coset of H in G does not in general divide (G, \*) into disjoint classes.*

*Proof:* We prove this only by the following example.
    Take $Q \setminus \{0\}$ to be a S-special definite group under multiplication; $H = Z^+ \subset Q \setminus \{0\}$ be a semigroup of $Q \setminus \{0\}$. Take $2 \in Q \setminus \{0\}$; $2Z^+ = \{2, 4, 6, \ldots\}$ is the S-definite special coset of the semigroup $Z^+$.

Take $\frac{1}{2} \in Q \setminus \{0\}$ and we find the S definite special coset of the semigroup $Z^+$.

$$\frac{1}{2} Z^+ = \left\{ \frac{1}{2},\ 1,\ \frac{3}{2},\ 2,\ \frac{5}{2},\ 3,\ \frac{7}{2}, \ldots \right\}.$$

Clearly $\frac{1}{2} Z^+ \cap 2Z^+ = 2Z^+$. So the S-definite special coset of the semigroup H of a group in general are not disjoint.

Now we proceed on to define the notion of Smarandache special definite double coset of a S-special definite group G.

**DEFINITION 2.1.10:** *Let G be a S-special definite group. Let H and K be two semigroups of G. For $x \in G$ we define Smarandache special definite double coset of the semigroups (S-special definite double coset of the semigroups) H and K for $x \in G$ as HxK = {hxk | $h \in H$ and $k \in K$}.*

We first illustrate this by an example before we prove and define any of the important properties.

*Example 2.1.16:* Let $Q \setminus \{0\}$ be a S-special definite group under multiplication. Let $K = 3Z^+$ and $H = 2Z^+$ be any two semigroups of $Q \setminus \{0\}$. To find the S-special definite double coset for $-1 \in Q \setminus \{0\}$. K – 1H = {–6, –12, –18, –24, –36, …}. We observe the following.

1.  $(K - 1 H) \cap H = \phi$;
    $(K - 1 H) \cap K = \phi$;
    But $K \cap H \neq \phi$ for {6, 12, 18, 24, …} are in $K \cap H$.
2.  K – 1 H is not even a closed set under multiplication.

Now we find the S-definite special double coset of H and K. K2H where $2 \in H \cdot K2H = \{12, 24, 36, \ldots\}$. The S-definite special double coset is a semigroup $12Z^+$ and in fact, $12Z^+$ is contained in both H and K as subsemigroup. Now we find the S-definite special double coset for $x = 5 \in Q^+ \setminus \{0\}$ using the semigroups H and K.

H5K = {30, 60, 120, 240, 360, …}. It is easily verified H5K is a subsemigroup of H and K.

As in case of S-special definite semigroup cosets we see S-definite special double cosets of a semigroup also need not divide the S-special definite group G into disjoint classes.

**THEOREM 2.1.7:** *Let G be a S-special definite group. H and K be semigroups of G. The S-special double definite coset of H and K need not be disjoint in general.*

*Proof:* We prove this by giving some examples. Let Q \ {0} be the group under multiplication. Take K = $2Z^+$ and H = $5Z^+$. Now take x = –1 ∈ Q\ {0} to find the S-special definite double cosets of the semigroups H and K for x = – 1. K – 1 H = – 10 $Z^+$ = {–10, –20, –30, –40, …}. Take –3 ∈ Q \ {0}. Find the S-definite special double coset of the semigroups H and K using –3.

K –3 H = {–30, –60, –90, …}. We see K -1 H ∩ K- 3H ≠ ϕ. Hence the claim. It is pertinent to mention here that we can have (H x K ) ∩ H y K = ϕ for some x and y. For take in the S-special definite group, H = 2Z and K = 3Z to be semigroups. Find H – 1 K for –1 ∈ Q \ {0} the S-special definite double coset of the semigroups H and K. H – 1 K = {–6, –12, –18, –24, …}. Now using 5 ∈ Q \ {0} we get the S-special definite double coset of the semigroups H and K to be H 5K = {30, 60, 90, 90, 120, …}. Clearly K –1 H ∩ H 5 K = ϕ. Hence the claim.

Thus unlike the double coset of a group the S-definite special double coset of the semigroups of a S-special definite group G do not partition G. This is shown from the examples just given in the theorem 2.1.7.

Another interesting factor about these S-special definite double cosets of semigroups of S special definite groups is that in many a cases they form a semigroup.

We prove an interesting theorem in this regard.

**THEOREM 2.1.8:** *Let Q \ {0} be the S-special definite group under multiplication. For any pair of semigroups of the form $mZ^+$ and $nZ^+$ of Q \ {0} and for all x ∈ $Z^+$ we see the S-definite* *special double coset of the semigroups $mZ^+$ and $nZ^+$ is a semigroup of $Q \setminus \{0\}$.*

*Proof:* Given $Q \setminus \{0\}$ to be a S-special definite group under multiplication. Let $mZ^+$ and $nZ^+$ be two semigroups of $Q \setminus \{0\}$. For and $x \in Z^+$. We consider $mZ^+ \times n Z^+ = mxnZ^+$ it is clear $(mxn) \in Z^+$ and $mxn Z^+$ is again a semigroup such that $mxn Z^+$ is contained in $mZ^+$ and $nZ^+$ i.e., $mxn Z^+$ is a subsemigroup of $mZ^+$ and $nZ^+$.

Further if $x \in Z^-$ or $x \in (Q \setminus \{0\}) \setminus Z^+$ then $m Z^+ x n Z^+$ is not even closed under the multiplication of rationals! This is evident for if $x \in Z^-$ then $m Z^+ x n Z^+ = mxnZ^+ =$ set of negative multiples of $mxn$ of the form $N(mxn)$ $N = \{1, 2, \ldots\}$ so t $mxn \cdot r m x n =$ t r m n would be + ve for any $t, r \in N$. Thus $mxn Z^+$ is not even closed under multiplication. Further if $x \in Q \setminus \{Z\}$ then also $m Z^+ x n Z^+ = (m \times n) Z^+$ is not closed under multiplication of rationals. Hence the claim.

However we illustrate this situation by the following example.

*Example 2.1.17:* Let $Q \setminus \{0\}$ be the group under multiplication. Let $5 Z^+$ and $6 Z^+$ be two semigroups of $Q \setminus \{0\}$. Take $8 \in Z^+$; $5Z^+ 8\ 6 Z^+ = 240 Z^+ = \{240, 480, 720, 960, \ldots\}$ is a semigroup. If we take $-3 \in Z^-$ and find $5Z^+(-3)\ 6Z^+ = \{-180, -360, -540, \ldots\}$. Now the product of any two elements in $5Z^+(-3)\ 6Z^+$ does not belong to $5Z^+(-3)\ 6Z^+ = \{-180, -360, -540, \ldots\}$. Now the product of any two elements in $5Z^+(-3)\ 6Z^+$ does not belong to $5Z^+(-3)\ 6Z^+$ as it gives a positive number. Now we take some $\frac{1}{7} \in Q \setminus \{0\}$. We find out

$$5Z^+(\frac{1}{7})\ 6Z^+ = \left(\frac{30}{7}, \frac{60}{7}, \frac{90}{7}, \ldots\right).$$

Now take

$$\left(\frac{30}{7} \times \frac{60}{7} = \frac{180}{49} \notin 5Z^+ \left(\frac{1}{7}\right) 6Z^+\right).$$

Thus $5Z^+ \frac{1}{7} 6Z^+$ is not a semigroup under product. Now we have seen the S-special definite double coset of S-definite special group.

We define the product of two semigroups of a S-special definite group.

**DEFINITION 2.1.11:** *Let (G, *) be a S-special definite group. H and K be any two semigroups of G. HK = {h * k / h ∈ H and k ∈ K} is defined as the product of H and K called the Smarandache definite special product of semigroups (S-definite special product of semigroups) in a S-special definite group (G, *).*

We propose a few problems to this end in the last chapter of this book. However we illustrate this concept by the following example.

*Example 2.1.18:* Let Q \ {0} be a S-definite special group under multiplication. Take H = $5Z^+$ and K = $3Z^+$; both H and K are semigroups of the S-definite special group Q \ {0}. Now HK = {15, 30, 45, 60, 75, …}. We see in this case HK is also a semigroup under multiplication.

## 2.2 Smarandache special definite algebraic structures

In this section for the first time we define a new Smarandache algebraic structure called Smarandache special definite fields and give some of its properties like their substructures.

**DEFINITION 2.2.1:** *Let (F, +, *) be a field. We call F a Smarandache special definite field (S-special definite field) if F contains a proper subset T such that (T, +, *) is just a ring.*

We give an illustrative example.

*Example 2.2.1:* Let Q be the field of rationals, Q is a S-special definite field as it contains Z to be a proper subset which is a ring.

All fields are not in general S-special definite fields.

**THEOREM 2.2.1:** *All fields are not in general S-special definite fields.*

*Proof:* We show there are fields which do not contain a proper P subset which is a ring. Take $Z_2 = \{0, 1\}$. $Z_2$ is not a S-special definite field. $Z_p = \{0, 1, 2, \ldots, p-1\}$ is a prime field of characteristic p (p a prime) has no proper subset P such that P is a ring. Hence the claim.

Let us define the Smarandache special definite subfield of a field G.

**DEFINITION 2.2.2:** *Let F be a field. Suppose K be a proper subfield of the field F and if K has a proper subset S ($\neq \phi$) such that S is a ring then we call K to be a Smarandache special definite subfield (S-special definite subfield) of the field F.*

We first give an example of the S-special definite subfield of the field F.

*Example 2.2.2:* Let R be the field of reals. $Q \subseteq R$ is the subfield of R and we see Q is a S-special definite subfield of R. For take T = 3Z; 3Z is a ring contained in Q.

In view of this we have the following theorem.

**THEOREM 2.2.2:** *Let F be a field; suppose F has a S-definite special subfield then F is also a S-definite special field.*

*Proof:* Given F is a field such that F has a subfield say K; where K is a S-definite special subfield of F. Thus K has a proper non empty set T such that T is a ring under the operations of K. Now $T \subset K$ but $K \subset F$ so $T \subset F$ and T is a ring so F is a S-definite special field.

Now we prove even if F be a S-definite special field then in general F need not have S-definite special subfield.

**THEOREM 2.2.3:** *Let F be a S-definite special field. F in general need not contain a S-definite special subfield.*

*Proof:* We prove this only by example 2.2.3.

*Example 2.2.3:* Take Q to be a field. Clearly Q is a S-definite special field. Now we know as Q is the prime field of characteristic zero; Q has no proper subfields so the question of Q having a subfield does not rise so Q having a S-definite special subfield is an impossibility.

Based on this we define the notion of Smarandache special definite prime field.

**DEFINITION 2.2.3:** *Let F be any field. If F has no S-special definite subfield then we call F to be a Smarandache special definite prime field (S-special definite prime field).*

*Example 2.2.4:* R the field of reals is not a S-special definite prime field.

*Example 2.2.5:* Let Q be field of rationals. Q is a S-definite special prime field.

*Example 2.2.6:* Let Q[x] be the ring of polynomials. Let $p(x) = x^2 + 1 \in Q[x]$. Now

$$\frac{Q[x]}{\langle p(x) \rangle} = T = \{ax + b \mid b, a \in Q\}$$

where $\langle p(x) \rangle$ is the ideal generated by $x^2 + 1$. $Q \subseteq T$; Q is a subfield of T and Q is a S-definite special subfield. Hence T is a S-definite special field.

We have the following interesting theorem.

**THEOREM 2.2.4:** *Let F be a prime field of characteristic zero then F has no S-definite special subfields.*

*Proof:* If F is a prime field, F has no subfields hence F cannot have S-definite special subfield.

We have the converse also to be true.

**THEOREM 2.2.5:** *Let F be a field of characteristic zero. If F has no S-definite special subfields then F is a prime field.*

*Proof:* We know if F is the field of characteristic zero then it contains properly or otherwise the field of rationals. Given F is a field of characteristic zero. We know if F has a subfield of characteristic zero say T then T can be Q or T is any field containing Q. So that T is a S-special definite subfield of F. Hence F is not prime. So if F is a prime field then F has no S-definite special subfields.

It is left as a open problem. Can a finite field be a S-definite special field?

*Remark:* It is important to note that any S-definite special field can have only rings which are always integral domains.

Now we proceed on to define the notion of S-definite special ideal related with the ring contained in S-definite special field.

**DEFINITION 2.2.4:** *Let F be a S-definite special field. Let R be a ring in F where R is not a field. An ideal I in R will be called as the Smarandache definite special ideal (S-definite special ideal) of F. We know a field can never have ideals but S-definite special fields can have S-definite special ideals.*

We first illustrate the following example.

*Example 2.2.7:* Let Q be the S-definite special field. Z is the ring contained in Q. pZ is a S-definite special ideal of Q. (p any positive number in N).

***Example 2.2.8:*** Let Q[x] be the polynomial ring. Let $p(x) = x^2 - 2$ and I the ideal generated by $x^2 - 2$. The quotient ring

$$\frac{Q[x]}{\langle x^2 - 2 \rangle = I} = V = \{ax + b + I \mid a, b \in Q\}$$

is a S-special definite field. Now take the ring $T = \{ax + b + I \mid a, b \in Z\}$; $T \subseteq V$. So V is a S-special definite field. If we take $W = \{ax + b + I \mid a, b \in 2Z\}$ then W is a S-definite special ideal of V.

Now we have seen examples of S-definite special ideals of a field F. Since the rings in a S-definite special fields are commutative all rings are commutative. So the question of right or left ideals does not arise. But we can have the notion of S-definite special maximal ideals, S-definite special minimal ideals, S-definite special prime ideals and S-definite special principal ideals.

**DEFINITION 2.2.5:** *Let F be a S-definite special field having V to be a subring in F. Let $P \subset V$ be an ideal of V if P is a maximal ideal then we call P to be a Smarandache definite special maximal ideal (S-definite special maximal ideal) of F if $P \subseteq M \subseteq V$, where M is an ideal of V containing P then $P = M$ or $M = V$. Likewise we define $T \subseteq V$ to be a minimal ideal of V as a Smarandache definite special minimal ideal (S-definite special minimal ideal) of F if $(0) \subseteq W \subseteq T$, W an ideal of V then $(0) = W$ or $W = T$.*

*Now we say a Smarandache definite special ideal P to be a Smarandache definite special principal ideal (S-definite special principal ideal) of F if P is generated by a single element.*

*Let M be a S-definite special ideal of F, if for $xy \in M$ then either $x \in M$ or $y \in M$ then we call M to be a Smarandache special definite prime ideal (S-special definite prime ideal) of F.*

Now we illustrate these by the following examples.

***Example 2.2.9:*** Let Q be the S-definite special field. Z be the ring contained in Q. Take $P = 6Z$ then P is a S-definite special ideal of Q which is not a S-definite special non prime ideal of Q as 2.3 ∈ P but 2 ∉ P and 3 ∉ P. Thus we have shown in general every S-definite special ideal of Q need not be a S-definite special prime ideal of Q.

Now consider T = 3Z then T is a S-definite special ideal of Q. Clearly 3Z is a S-definite special maximal ideal of Q. Consider the S-definite special ideal W = 2Z of Q then W is a S-definite special principal ideal of Q as W is generated by 2.

Clearly Q has no S-definite special minimal ideals.

We propose the problem of finding S-definite special minimal ideals of a field F.

Now we have defined the concept of S-definite special minimal ideals, S-definite special maximal ideals, S-definite special prime ideals and S-definite special principal ideals and exhibited them with examples.

We now proceed on to define notions analogues to galois group of automorphisms for S-definite special field.

**DEFINITION 2.2.6:** *Let F be a field suppose F(b) be the finite extension of F. Let V be the subring of F and V(b) be the finite extension subring of V.*

*Now G(F(b): F) is the group of automorphisms of F(b) fixing F. Let $G_s(V(b), V)$ be the group of automorphism of V(b) fixing V. We call $G_s(V(b), V)$ to be the Smarandache group of definite special automorphism of V(b) (S-group of definite special automorphism of V(b)) fixing V.*

We propose the following, Is $G_s(V(b), V) \cong G(F(b): F)$? Now we have only written in terms of abstract extensions. We would now define more in an non abstract level. We know the element a ∈ K is said to be algebraic of degree n over F if it satisfies a non zero polynomial of lower degree. If a ∈ K is algebraic of degree n over F then [F(a): F] = n where F(a) is {$\beta_0 + \beta_1 a + \ldots + \beta_{n-1} a^{n-1} \mid \beta_0, \beta_1, \ldots, \beta_{n-1} \in F$}.

Suppose R is in the ring in F such that $\alpha_0 a^n + \alpha_1 a^{n-1} + \ldots + \alpha_{n-1} a + \alpha_n = 0$; $\alpha_0, \ldots, \alpha_n \in R$ then R(a) = {$\beta_0 + \beta_1 a + \ldots + \beta_{n-1} a^{n-1} \mid \beta_0, \beta_1, \ldots, \beta_{n-1} \in R$}. R(a) ⊂ F(a) i.e., R(a) is a ring in F(a). R(a) is a module over the ring R. Now just analogous to [F(a): F ] we call [R(a): R], Smarandache definite special finite extension (S-definite special finite extension).

Now for every group of automorphisms of F(a) which keeps F fixed denoted by G(F(a); F); we define $G_s(R(a), R)$ to be the Smarandache definite special group of automorphisms (S-definite special group of automorphisms) of R(a) keeping R fixed.

Find the relation between $G_s(R(a); R)$ and G(F(a), F ). Several interesting results in this direction can be determined; in order to make these concepts more easy we would be defining the notion of Smarandache polynomial rings.

**DEFINITION 2.2.7:** *Let F be a field of characteristic zero. F(x) be the smallest field containing F and x. x a variable i.e., x is algebraic over F. Let Z be the ring in F. Then we have Z[x], the polynomial ring with coefficients in Z. Clearly Z[x] is a ring contained in F(x). We call Z[x] the Smarandache algebraic polynomial ring of the extension field (S-algebraic polynomial ring of the extension field) F(x) of F.*

*It is noted that this variate x is algebraic over F of degree some m (m > 1), i.e., we have $\alpha_0, ..., \alpha_m \in F$ such that $\alpha_0 x_m + ... + \alpha_m = 0$.*

We illustrate this situation by an example so that the reader may not have difficulty in understanding this concept.

***Example 2.2.10:*** Let Q be the field of rationals; $\sqrt{3}$ is algebraic over Q for it splits in the field $Q(\sqrt{3})$ and $Q(\sqrt{3})$ is the splitting field of the polynomial $p(x) = x^2 - 3$. Or more in the conventional language we can say $\sqrt{3} \in R$ (R – the real field) is algebraic over Q for $1.(\sqrt{3})^2 - 3 = 0$ and $[Q(\sqrt{3}): Q] = 2$, degree of the lowest degree of the polynomial which is satisfied by $\sqrt{3}$.

We say $Z(\sqrt{3})$ is a subring of $Q(\sqrt{3})$ and $Z$ is the subring of $Q$. $Z(\sqrt{3})$ is the Smarandache algebraic polynomial ring of the extension field $Q(\sqrt{3})$.

In this way we can construct any number of examples. Also $Z(\sqrt{3})$ is the module over $Z$.

We shall be defining these concepts when we relate it with a field as Smarandache concepts in the following chapter.

We now proceed on to define Smarandache special definite homomorphism.

**DEFINITION 2.2.8:** *Let F and K be any two S-definite special fields having R and T as subrings of F and K respectively. The Smarandache special homomorphism of S-special definite field (S-special homomorphism of S-special definite field) $\phi$ is a map from $R \subset F \to T \subset K$ such that $\phi$ is a ring homomorphism from R to T i.e.,*

$$\phi(r \pm s) = \phi(r) \pm \phi(s)$$
$$\phi(r\,s) = \phi(r)\,\phi(s)$$

*for all $r, s \in R$.*

*Note:* We are not bothered about $F \setminus R$ and $K \setminus T$ i.e., $\phi$ is not even defined on the whole of F only on the ring R of F.

We illustrate this by the following example.

*Example 2.2.11:* Let Q be the field of rationals and R the field of reals. Clearly Q and R are S-special definite fields. Let 2Z be the ring contained in Q and Z the ring contained in R. A map $\phi$ from $2Z \subset Q$ to $Z \subset R$ be defined by $\phi(n) \mapsto n$ for $n \in 2Z$ as $2Z \in Z$; clearly $\phi$ is a S-definite special homomorphism of S-definite fields Q to R, the field of reals.

Now one can define not only for any ring R contained in the field F to any other ring T contained in K but also define many homomorphisms from R to T itself. The notion of automorphism, isomorphism follow as in case of general ring homomorphism but in this case restricted only to the rings in them.

We illustrate this yet by another example.

***Example 2.2.12:*** Let $Q(\sqrt{2}, \sqrt{3})$ be the S-special definite field and Q be another S-special definite field. Let $Z[\sqrt{2}]$ be the ring in $Q(\sqrt{2}, \sqrt{3})$ and $2Z[\sqrt{2}]$ be the ring in $Q(\sqrt{2}, \sqrt{7})$.

The reader is expected to define a S-special definite homomorphism from $Q(\sqrt{2}, \sqrt{3})$ to $Q(\sqrt{2}, \sqrt{7})$. Clearly kernel $\phi$ is an ideal in the ring which we define as a S-special definite ideal of the S-special definite field F.

Now we have the notion of Smarandache special definite division rings.

**DEFINITION 2.2.9:** *Let K be a division ring. We say K is a Smarandache special definite division ring (S-special definite division ring) if K contains a nontrivial ring R which is not a division ring.*

We give the following example.

***Example 2.2.13:*** Let K be the division ring of real quaternions where F is the field of reals and $K = \{\alpha_o + \alpha_1 i + \alpha_2 j + \alpha_3 k \mid \alpha_0, \alpha_1, \alpha_2, \alpha_3 \in F$ and $i^2 = j^2 = k^2 = -1$, $ij = k$, $jk = I$, $ki = j$, $ji = -k$, $kj = -i$ and $ik = -j\}$.
Define multiplication in K as follows.

$$(\alpha_0 + \alpha_1 i + \alpha_2 j + \alpha_3 k)(\beta_0 + \beta_1 i + \beta_2 j + \beta_3 k)$$
$$= \gamma_o + \gamma_1 i + \gamma_2 j + \gamma_3 k$$

where
$$\gamma_0 = \alpha_0 \beta_0 - \alpha_1 \beta_1 - \alpha_2 \beta_2 - \alpha_3 \beta_3$$
$$\gamma_1 = \alpha_0 \beta_1 + \alpha_1 \beta_0 + \alpha_2 \beta_3 - \alpha_3 \beta_2$$
$$\gamma_2 = \alpha_0 \beta_2 + \alpha_2 \beta_0 - \alpha_1 \beta_3 + \alpha_3 \beta_1$$

$$\gamma_3 = \alpha_0\beta_3 - \alpha_2\beta_1 + \alpha_1\beta_2 + \alpha_3\beta_0.$$

It is easily verified that K is a division ring. Take $Z \subset K$, Z is just a ring and not a division ring so K is a S-special definite division ring.

We can also take $T = \{\alpha_o + \alpha_1 i + \alpha_2 j + \alpha_3 k \mid \alpha_o, \alpha_1, \alpha_2, \alpha_3 \in Z\}$; T is a ring in K, which is a non commutative ring.

As in case of S-special definite fields we define in case of S-special definite division rings S-special definite homomorphisms of division rings.

**DEFINITION 2.2.10:** *Let K and T be any two S-special definite division rings. We call $\phi$ a Smarandache special definite division rings homomorphism (S-special definite division rings homomorphism) if $\phi$ a ring homomorphism from a ring R in K to a ring S in T.*

It is pertinent to mention here that the map $\phi$ will not be defined on the whole of K but only on the whole of the ring for which we need the ring homomorphism. As in case of ring homomorphism even in case of S-definite special definite division rings we have ker $\phi$ to be an ideal of the ring R in K. We have proposed some problems of finding the S-definite special homomorphisms of S-definite special division ring.

Now having defined the notion of S-special definite fields and S-special definite division rings we now in the next chapter define S-special definite modules using S-special division rings and also Smarandache definite special vector spaces using S-special definite fields. We will also define the notion of S-special definite linear algebra.

**DEFINITION 2.2.11:** *Let (G, \*) and (H, o) be any two S-special definite groups. We define $\phi$ to be a Smarandache special definite group homomorphism (S-special definite group homomorphism) of the S-special definite groups if $\phi$ is a semigroup homomorphism from $S_1$ to $S_2$ where $S_1$ is a semigroup in G and $S_2$ is a semigroup in H; i.e., $\phi(x * y) = x \circ y$ for all $x, y \in S$.*

We illustrate this situation by the following example.

***Example 2.2.14:*** Let $Q \setminus \{0\}$ be the S-special definite group and

$$M_{2 \times 2} = \left\{ \begin{pmatrix} a & b \\ c & d \end{pmatrix} \middle| a, b, c, d \in Q \text{ and } ad - bc \neq 0 \right\}$$

be the S-special definite group. Let 2Z be the semigroup under multiplication and

$$T = \left\langle \begin{pmatrix} 1 & -2 \\ -2 & 1 \end{pmatrix} \right\rangle$$

be the semigroup in $M_{2 \times 2}$ generated by

$$\begin{pmatrix} 1 & -2 \\ -2 & 1 \end{pmatrix}.$$

Define a non trivial semigroup homomorphism from 2Z to T.

**Chapter Three**

# SMARANDACHE SPECIAL DEFINITE RINGS AND SMARANDACHE DEFINITE SPECIAL VECTOR SPACES

This chapter has 3 sections. In section one the new notion of S-special definite rings are introduced. Section two introduces the notion of S-special definite vector spaces. It is shown the notion of Gram Schmidt orthogonalization process cannot be defined on these S-special vector spaces. In the final section the notion of S-special definite linear algebras are defined and some interesting properties derived.

## 3.1 Smarandache Special Definite Ring

In this section we for the first time define the notion of Smarandache special definite rings and give some of its interesting properties.

**DEFINITION 3.1.1:** *Let R be a ring, we say R is a Smarandache special definite ring (S-special definite ring) if R has a proper subset $S \neq 0$ such that S is a semiring.*

*If S is a commutative semiring then we say R to be a Smarandache special definite weakly commutative ring(S-special definite weakly commutative ring). Only when every S is a commutative semiring then alone we call R to be a S-special definite commutative ring.*

We first give some examples of them.

***Example 3.1.1:*** Let Z be the ring of integers. Take $S = 2Z^+ \cup \{0\}$ a, semiring; clearly $S \subset Z$, so Z is a S-special definite ring. Since Z happens to be commutative we can say Z is trivially a S-special definite commutative ring.

***Example 3.1.2:*** Let $M_{2\times 2}$ be the collection of all $2 \times 2$ matrices with entries from Q. $M_{2\times 2}$ is a ring under matrix addition and matrix multiplication.
Take
$$S = \left\{ \begin{pmatrix} a & b \\ c & d \end{pmatrix} \middle| a,b,c,d \in Z^+ \cup \{0\} \right\}$$

a proper subset in $M_{2\times 2}$. S is a semiring with respect to matrix addition and multiplication. But S is not a commutative semiring. So $M_{2\times 2}$ is only a S-definite special non commutative ring.
Take
$$P = \left\{ \begin{pmatrix} a & 0 \\ 0 & 0 \end{pmatrix} \middle| a \in 2Z^+ \cup \{0\} \right\}$$

a proper subset in $M_{2\times 2}$. P is a commutative semiring, hence $M_{2\times 2}$ is a S-definite special weakly commutative ring. For $M_{2\times 2}$ has both commutative semiring as well as non commutative semiring so $M_{2\times 2}$ can only be a S-special definite weakly commutative ring.

We have the following interesting theorem.

**THEOREM 3.1.1:** *Every commutative ring R is a S-special definite commutative ring.*

*Proof:* Since every subset S of R which is a semiring happens to be commutative we see that R is a S-special definite commutative ring.

**DEFINITION 3.1.2:** *Let R be a ring we say a subring V of R to be a Smarandache special definite subring (S-special definite subring) of R if V has a peroper subset $W \neq \phi$ such that W is a semiring.*

We have the following interesting theorem.

**THEOREM 3.1.2:** *Let R be a ring. If R has a Smarandache special definite subring then R itself is a S-special definite ring.*

*Proof:* Given R is a ring such that R contains a proper subring V such that V is a S-special definite subring of R. Now given V is a S-special definite subring so V contains a nonempty subset W such that W is a semiring. Now $W \subset V$ and W is a semiring; but V is a subring of R so $V \subset R$ and $W \subset V \subset R$. Thus $W \subset R$ is a semiring. Hence R is a S-definite special ring.

*Example 3.1.3:* Let R be the ring of reals. Q is a subring of R. Now Q contains a proper subset $Z^+ \cup \{0\} = Z^\circ$ such that $Z^\circ$ is a semiring. Hence Q is a S-special definite subring of R. Clearly as $Z^\circ \subset Q \subset R$ and $Z^\circ$ is a semiring we see R is also a S-special definite subring.

*Example 3.1.4:* Consider Z to be the ring 2Z is a subring of Z. 2Z is a S-special definite subring of Z as 2Z contains a proper subset $V = 2Z^+ \cup \{0\}$ which is a semiring. Thus Z is itself a S-special definite ring.

Now we proceed, on to define the notion of Smarandache special definite ideal of a ring.

**DEFINITION 3.1.3:** *Let R be any ring we call a proper subset I of R to be a Smarandache special definite ideal (S-special definite ideal) of R if the following conditions are satisfied.*
  *(1) I is a S-special definite subring of R.*
  *(2) I is an ideal of R.*

We illustrate this situation by the following examples.

*Example 3.1.5:* Let Z be the ring of integers. Take I = 5Z, I is a S-special definite ideal of Z as 5Z is a S-special subring of R and V = $5Z^+ \cup \{0\}$ is a proper subset of 5Z which is a semiring of Z.

*Example 3.1.6:* Let

$$M_{2\times 2} = \left\{ \begin{pmatrix} a & b \\ c & d \end{pmatrix} \middle| a,b,c,d \in Z \right\}$$

be the collection of all 2 × 2 matrices with entries from Z, $M_{2\times 2}$ is a S-definite special ring. Further if we take

$$I_{2\times 2} = \left\{ \begin{pmatrix} a & b \\ c & d \end{pmatrix} \middle| a,b,c,d \in 2Z \cup \{0\} \right\},$$

$I_{2\times 2} \subseteq M_{2\times 2}$ and $I_{2\times 2}$ is a S-special definite ideal of $M_{2\times 2}$. For

$$V = \left\{ \begin{pmatrix} a & b \\ c & d \end{pmatrix} \middle| a,b,c,d \in 2Z^+ \cup \{0\} \right\}$$

is a semiring contained in $I_{2\times 2}$. Hence the claim.

It is interesting to see the following the results. As the proof of them are straight forward we just only mention the results.

*Result I:* Every S-definite special ideal of a ring R is a S-definite special subring of R.

This is by the very definition of S-definite special ideal of a ring R.

*Result II:* If a ring R has a S-definite special ideal then the ring R is a S-definite special ring.

This proof is also straight forward using result I the result follows.
    Now we have certain rings which can have S-definite special subrings which are not S-special definite ideals of R.

We illustrate this situation by the following example.

*Example 3.1.7:* Let Q be the ring of rationals. Z is a subring of Q. Clearly Z is a S-special definite subring of Q but is not a S-special definite ideal of Q. Hence the claim.

Now we illustrate a special ring in which
(1) Every subring is a S-special definite subring and
(2) Every S-definite special subring is also a S-definite special ideal of R.

**THEOREM 3.1.3:** *Let Z be the ring of integers. Z is a S-definite special ring such that*
    *(1) Every S-special definite subring of Z is a S-special definite ideal of Z.*
    *(2) Every subring of Z is a S-definite special subring.*

*Proof:* Given Z is a ring, clearly Z is a principal ideal domain. Every ideal is principal. The subrings of Z are only of the form nZ where n is a positive integer. Thus every subring is an ideal. Every ideal in Z is only of the form nZ as Z is a principal ideal domain. So every subring of Z is a S-special subring of Z and clearly all S-special subrings of Z are S-special definite ideals of Z. It is easily verified every subring in Z is generated only by a single element. For if S is a proper subring of Z and by chance S is not generated by a single element, so that S is generated by atleast two elements say such that one of them is the first least and another the second least, let it be x and y (without loss of generality we can assume S is generated by at least two elements). Further we can assume x < y. Now x, y ∈ S so x + y and x – y ∈ S. But now if x – y = t > 0 then t ≤ x but x is the least element of S which gives a contradiction if t < x. If t = x then we see y = 2x i.e. y = x + x so x is the only generator. Thus every subring of Z is generated by a single element and none of the subrings S, contain 1 for if it contains 1, we have S = Z i.e. S is not a proper subring of Z.

Every subring is obviously a S-special definite subring of Z for if S is subring it is only of the form nZ (n > 1) now take $V = nZ^+ \cup \{0\}$, V is a semiring; hence S is a S-definite special subring of Z.

We define in view of this property Smarandache strong special definite rings.

**DEFINITION 3.1.4:** *Let R be any ring we call R to be Smarandache strong special definite ring (S-strong special definite ring) if every subring of R is a S-special definite subring of R.*

*Example 3.1.8:* Let Z be the ring of integers; Z is a S-strong special definite ring.

**DEFINITION 3.1.5:** *Let R be a ring we call R to be a Smarandache ideally strong definite special ring (S-ideally strong definite special ring) if every S-definite special subring is a S-definite special ideal of R.*

Cleary Z the ring of integers is a S-ideally strong definite special ring.
   We see all group rings over the field of characteristic zero are S-definite special rings.

For literature about group rings refer [70-1].

**THEOREM 3.1.4:** *Let F be a field of characteristic zero and G any group. The group ring FG is a S-definite special ring.*

*Proof*: Given F is a field of characteristic zero so either F = Q or Q ⊆ F in both cases we see $Z^+ \cup \{0\}$ is a proper subset of F. Now FG is the group ring of G over F so F ⊆ FG (∴ 1 ∈ G). Also V = $Z^+ \cup \{0\}$ ⊂ F ⊂ FG. Thus FG contains a proper subset V such that V is a semiring.

Hence FG is a S-definite special ring when F is of characteristic zero.

*Example 3.1.9:* Take Q to be the field of rationals G = $S_3$, the symmetric group of degree 3. QG is a S-definite special ring for V = $Z^+ \cup \{0\}$ is a semiring also ($Z^+ \cup \{0\}$)G = VG, the group semiring is also a semiring.

**THEOREM 3.1.5:** *All group ring ZG of the group G over the ring of integers Z is a S-definite special ring.*

*Proof:* Since V = ($Z^+ \cup \{0\}$) ⊆ ZG we see ZG is a S-special definite rings. Clearly V is a semiring. It has also infinite number of semirings given by ($pZ^+ \cup \{0\}$).

The group ring ZG has infinite number of S-definite special ideals.

**THEOREM 3.1.6:** *Let ZG be the group ring of the group G over the ring of integers Z. ZG has S-special definite ideal.*

*Proof:* ZG is the group ring; $2Z^+ \cup \{0\}$ = V ⊆ Z ⊆ ZG. V is a semiring, consider ($2Z^+ \cup \{0\}$)G = VG, VG is an ideal of ZG. Hence the claim.
We have given a nontrivial class of S-definite special rings.

Now we give yet another class of S-definite special rings.

Let K be any field or a commutative ring with unit. S any semigroup with unit or a monoid. The semigroup ring KS of the semigroup over K (the field or ring) is defined analogous to group rings. KS the semigroup ring is also a ring. For more about semigroup rings please refer [101]. We give some special properties about semigroup ring. If K is a field of characteristic zero and S any semigroup with unit, then the semigroup ring is S-special definite ring which is proved in the form of the following theorem.

**THEOREM 3.1.7:** *Let K be a field of characteristic zero. S be any semigroup with 1. Then the semigroup ring KS is a S-special definite ring.*

*Proof:* Let KS be the semigroup ring where K is a field of characteristic zero and S is a semigroup with 1. To prove KS is a S-special definite ring; i.e., KS has a proper subset V where V is a semigroup. Since we are given K is a field of characteristic zero either Q = K or Q ⊆ K, according as K is a prime field of characteristic zero are a non prime field. Clearly Z ⊂ Q ⊆ K; now take V = $Z^+$ ∪ {0}, a proper subset of Z. V is a semiring; also VS is a semiring. Hence KS is a S-sepcial definite ring.

Thus we have given yet another class of S-special definite rings. Throughout this book by S(n) we denote the semigroup of all mappings of the set (1 2 3 … n) to itself under the composition of mappings. How ever we express this S(n) by some examples.

*Example 3.1.10:* Let S(2) be the semigroup of mappings of (12) to itself.

$$S(2) = \left\{ \begin{pmatrix} 1 & 2 \\ 1 & 2 \end{pmatrix} = e, \begin{pmatrix} 1 & 2 \\ 2 & 1 \end{pmatrix} = p_1, \begin{pmatrix} 1 & 2 \\ 1 & 1 \end{pmatrix} = p_2, \begin{pmatrix} 1 & 2 \\ 2 & 2 \end{pmatrix} = p_3 \right\};$$

S(2) is a semigroup under the composition of maps for

$$p_1 : \begin{matrix} 1 \to 2 \\ 2 \to 1 \end{matrix},$$

$$p_2 : \begin{matrix} 1 \to 1 \\ 2 \to 1 \end{matrix},$$

$$p_3 : \begin{matrix} 1 \to 2 \\ 2 \to 2 \end{matrix}$$

and

$$e: \begin{matrix} 1 \to 1 \\ 2 \to 2 \end{matrix};$$

i.e., the identity map $p_3 \circ e = e \circ p_3 = p_3$,

$$p_3 \circ p_3 = \begin{matrix} 1 \to 2 \\ 2 \to 2 \end{matrix} = p_3,$$

$$p_2 \circ p_2 = \begin{matrix} 1 \to 1 \\ 2 \to 1 \end{matrix} = p_2,$$

$$p_1 \circ p_2 = \begin{matrix} 1 \to 1 \\ 2 \to 1 \end{matrix}$$

and so on. $S(2)$ is a semigroup of order $2^2$ usually called as symmetric semigroups.

Please refer [101] for more information.

**Example 3.1.11:** Take $S(4)$ to be the symmetric semigroup. $S(4)$ to be the symmetric semigroup. $S(4)$ has $4^4$ elements in it.

Thus we see $S(n)$ gives us the symmetric semigroup of finite order when $n < \infty$ and the number of elements in $S(n)$ is $n^n$.

**Example 3.1.12:** Take $Z_2 = \{0, 1\}$, the prime field of characteristic two and $S(2) = \{e, p_1, p_2, p_3\}$ be the symmetric semigroup got using (1, 2). Now $Z_2 S(2)$ be the semigroup ring. $Z_2 S(2) = \{0, 1 = e, p_1, p_2, p_3, 1 + p_1, 1 + p_2, 1 + p_3, p_1 + p_2, p_1 + p_3, 1 + p_1 + p_2, 1 + p_1 + p_3, 1 + p_2 + p_3, p_1 + p_2 + p_3, 1 + p_1 + p_2 + p_3\}$ Clearly $Z_2 S(2)$ is not a S-definite special ring.

**Example 3.1.13:** Let $ZS(2)$ be the semigroup ring. $ZS(2)$ is a S-definite special ring.

Now having seen a few classes of S-definite special rings we now proceed on to define some more properties about S-definite special rings.

We give yet another class of non commutative rings which are S-special definite rings.

**THEOREM 3.1.8:** *Let $M_{n \times n} = \{m = (m_{ij}) | m_{ij} \in Z\}$; i.e., the set of all $n \times n$ matrices with entries form Z. $M_{n \times n}$ is a ring under matrix addition and matrix multiplication. $M_{n \times n}$ is a S-definite special ring.*

*Proof:* To show $M_{n \times n}$ is a S-definite special ring we must find a non empty subset V in $P_{n \times n}$ which is only a semiring. Take $V = P_{n \times n} = \{p = (p_{ij}) | p_{ij} \in Z^+ \cup \{0\}\}$, $P_{n \times n}$ is a semiring clearly V is a proper subset of $M_{n \times n}$; hence $M_{n \times n}$ is a S-definite special ring.

**THEOREM 3.1.9:** *Let $F_{n \times n}$ denote the set of all $n \times n$ matrices with entries from a field F of characteristic zero. $F_{n \times n}$ is a S-definite special ring.*

*Proof:* $F_{n \times n} = \{A = (a_{ij}) | a_{ij} \in F\}$ where F is a field of characteristic zero}. $F_{n \times n}$ is a ring under matrix addition and matrix multiplication. Take $V = \{V = (v_{ij}) | v_{ij} \in Z^+ \cup \{0\}\}$. Clearly $Z^+ \cup \{0\}$ is a proper subset of F as F is given to be a field of characteristic zero so F contains either Q or F = Q itself. $V = V_{n \times n}$ is a proper subset of $F_{n \times n}$ and $V = V_{n \times n}$ is a semiring. Hence $F_{n \times n}$ is a S-definite special ring. Thus we have seen yet another class of S-definite special rings.

*Example 3.1.14:* Let $M_{3 \times 3} = \{M = (m_{ij}) | m_{ij} \in R$, the field of reals $\}$, be the ring. Take $T_{3 \times 3} = \{T = (t_{ij}) | t_{ij} \in Z\}$; $T_{3 \times 3}$ is only a S-definite special subring of $M_{3 \times 3}$ which is never an ideal of $M_{3 \times 3}$.

We have yet another large class of S-definite special rings.

**THEOREM 3.1.10:** *Let F be a field of characteristic zero; F[x] the polynomial ring in the variable x. F[x] is a S-definite special ring.*

*Proof:* F[x] is a polynomial ring and $F \subseteq F[x]$. Also F is given to be field of characteristic zero so F either contains Q, the field of rationals or F = Q. In either case we have $Z^+ \cup \{0\} = V$, is a proper subset of F hence a proper subset of F[x]. But both V and V[x] are semirings hence F[x] is a S-definite special ring. Apart from this Z[x] is also a S-definite special ring.

Now we show that Z[x] can have S-special definite subrings which are not S-special definite subrings and which are not S-special definite ideals of Z[x].

*Example 3.1.15:* Let Z[x] be the polynomial ring over Z. Let

$$T[x] = \left\{ \sum_{i=D}^{n} a^i x^{2i} \,\middle|\, a_i \in Z \right\}.$$

T[x] is a subring of Z[x]. Clearly T[x] is a S-definite special subring of Z[x] as $Z \subseteq T[x]$ and Z in turn contains $V = Z^+ \cup \{0\}$ and V is a semiring. We see T[x] is not even an ideal of Z[x] as x. $T[x] \not\subseteq T[x]$, hence the claim!

Now we give yet another class of S-definite special rings. Let $K[x_1, x_2, \ldots, x_n]$ be a polynomial ring in the n variables $x_1, x_2, \ldots, x_n$ where K is a field or a commutative ring with 1.

For more about these rings refer I.N.Heristen [35-8]. We show $K[x_1, x_2, \ldots, x_n]$ is a S-definite special ring if K is a field of characteristic zero.

**THEOREM 3.1.11:** *Let K be a field of characteristic zero. $K[x_1, x_2, \ldots, x_n]$ be a polynomial ring in the n-variables. $K[x_1, x_2, \ldots, x_n]$ is a S-definite special ring.*

*Proof:* Given $K[x_1, x_2, \ldots, x_n]$ is a polynomial in the n-variables with coefficients from K, K a field of characteristic zero. Since K is a field of characteristic zero we have $Q = K$ or $Q \subseteq K$. Hence $V = Z^+ \cup \{0\} \subset K$ is a semiring in K thus $V \subseteq K[x_1, \ldots, x_n]$ is a S-definite special ring.

In fact $V[x_1, \ldots, x_n]$ is also a semiring in $K[x_1, x_2, \ldots, x_n]$ we have the following interesting theorem.

**THEOREM 3.1.12:** *Let R be a ring which is a S-definite special ring. Then R[x] is also a S-definite special ring.*

*Proof:* Given R is a S-definite special ring, so R contains a proper subset P such that P is a semiring. Since R ⊂ R[x] we have P ⊂ R[x] and P is a semiring in R[x] so R[x] a S-definite special ring.

However if R[x] be any polynomial ring over a ring R. If R[x] is a S-definite special ring will R be a S-definite special ring?

It is important at this juncture to mention that in general that a ring which is a S-definite special ring need not be a S-ring. Also all S-rings need not be a S-definite special rings but however there are rings which are both S-definite special rings as well as S-rings.

We now illustrate this situation by some theorems.

**THEOREM 3.1.13:** *Let Z[x] be a polynomial ring over the ring of integers Z. Z[x] is a S-definite special ring but Z[x] is not a S-ring.*

The proof is left as an exercise for the reader. For literature about S-rings please refer [101].

The following theorems are also left for the reader to prove as the proofs are straight forward.

**THEOREM 3.1.14:** *Let $S_{n \times n} = \{m = (m_{ij}) \mid m_{ij} \in Z\}$ be the ring of $n \times n$ matrices with entries form Z. $S_{n \times n}$ is a S-definite special ring which is not a S-ring.*

**THEOREM 3.1.15:** *Let $Z[x_1, x_2, \ldots, x_n]$ be a polynomial ring with coefficients from the ring of integers Z in the n-variables $x_1, x_2, \ldots, x_n$. $Z[x_1, x_2, \ldots, x_n]$ is a S-definite special ring and is not a S-ring.*

**THEOREM 3.1.16:** *Let $Z_n = \{\bar{0}, \bar{1}, \ldots, \overline{n-1}\}$ be the ring of integers modulo n (n a composite number), $Z_n$ is a S-ring and is not a S-definite special ring.*

Thus we see we have a class of rings which are S-rings but not S-definite special rings and still we have a class of rings which are S-definite special rings but not S-rings.

Now we will give a class of rings which are both S-rings as well as S-definite special rings which will answer the question that we have a class of rings which is both a S-ring and S-definite special ring.

**THEOREM 3.1.17:** *Let Q[x] be a polynomial ring Q[x] is both a S-ring as well as a S-definite special ring.*

**THEOREM 3.1.18:** *Let $M_{p \times p} = \{M = (m_{ij}) \mid m_{ij} \in Q\}$ be the ring of $p \times p$ matrices with entries from the field of rationals. $M_{p \times p}$ is a S-ring and also $M_{p \times p}$ is a S-special definite ring.*

*Proof:* Given $M_{p \times p} = \{M = (m_{ij}) \mid m_{ij} \in Q\}$ is the ring of all $p \times p$ matrices with entries from Q. Clearly $M_{p \times p}$ is a ring under matrix addition and matrix multiplication.

$M_{p \times p}$ is a S-ring. For take all matrices $M^1 = (m_{11}^1)$; i.e., $m_{ij}^1 = 0$ for $i \neq j$, $1 \leq j \leq p$ only $m_{11}$ is non zero; $T = M_1 \cup (0) = \{(m_{11}^1) \mid m_{11}^1 \in Q\} \cup (0)$, where (0) denotes the zero $p \times p$ matrix. Clearly T is a field under matrix addition and matrix multiplication. Infact $T \cong Q$ Thus $T \subset M_{p \times p}$ and as T is a field $M_{p \times p}$ is a S-ring. Now consider the subset $P \neq \phi$ in $M_{p \times p}$; where $P = \{A = (a_{ij}) \mid a_{ij} \in Z^+ \cup \{0\}\}$. We see P is a semiring under matrix addition and matrix multiplication; hence $M_{p \times p}$ is a S-definite special ring. Thus we see $M_{p \times p}$ is a S-ring as well as a S-definite special ring.

**THEOREM 3.1.19:** *The ring of polynomials Q[x] is a S-ring as well as a S-definite special ring.*

*Proof:* Now Q[x] is the ring of polynomials in the variable x with coefficients from the field Q. Q[x] is a S-ring for $Q \subseteq Q[x]$ and Q is a field. Further Q[x] is a S-definite special ring as $Z^+ \cup \{0\} = V \subset Z \subset Q \subset Q[x]$. Thus V is a proper subset of Q[x] and V is a semiring, hence the claim. Thus Q[x] is a S-ring as well as a S-definite special ring.

Those rings which are both S-rings as well as S-definite special rings gives the researchers more opportunity to study their structures. We give a new name to this class.

**DEFINITION 3.1.6:** *Let R be a ring we call R to be Smarandache doubly strong (S-doubly strong) ring if R is a S-ring as well as R is a S-definite special ring.*

We illustrate it by an example.

*Example 3.1.16:* The ring $Q[x_1, x_2]$ is a S-doubly strong ring.

Now we have the following interesting theorem which can be viewed as the existence theorem.

**THEOREM 3.1.20:** *The class of S-doubly strong rings is non empty.*

*Proof:* To show the class of S-doubly strong rings is non empty it is sufficient if we prove a class of rings which are S-doubly strong rings. Consider $M_{n \times n} = \{M = (m_{ij}) \mid m_{ij} \in R\}$. $M_{n \times n}$ is a ring under matrix multiplication and matrix addition which can be easily verified to be a S-ring as well as a S-definite special ring.
    This class of rings is non commutative. $Q[x_1, x_2, \ldots, x_n]$ is a S-doubly strong rings which form a class of commutative rings.
    Thus we see the class of S-doubly strong rings can be commutative or non commutative; hence we have shown the existence of S-doubly strong rings which infact forms an infinite class of rings.

So given a ring we can say either they are S-doubly strong rings or not i.e. we can divide the class of rings into two disjoint classes; we cannot always say all infinite rings are S-doubly strong rings. For we see $Z[x_1, x_2, \ldots, x_n]$ is not a S-doubly strong ring but $Z[x_1, x_2, \ldots, x_n]$ is a ring of infinite cardinality.

Thus we have studied a new class of rings which has a weaker substructure embedded in it. Now we have defined those class of rings as Smarandache special definite rings. Unlike Smarandache rings which contain a stronger structure in them here a reverse study is made. It is surprising to see we have a class of rings which are both S-rings as well as S-definite special rings. Further we see when we seek a weaker structure in a strong structure, like finding semigroups as a proper subset in a group or a semiring in a ring; we found that those algebraic structures invariably had infinite cardinality. We were unable to find such structures with finite number of elements in them. Further we are not in a position to say why structures with finite cardinality fail to cater to these S-definite special substructures. Does it imply when finite sets satisfy any algebraic operation those sets become very tight or dense or complete under those operations? For we could not find any finite group which had a proper subset which was a semigroup. Likewise we were unable to find finite rings which had in them a proper subset which was a semiring. We have left this as an open problems for the reader to prove or disprove this claim! How ever we feel it is impossible because of Cauchy theorems for groups.

Now we proceed on to define yet another Smarandache definite special algebraic structure for fields.

**DEFINITION 3.1.7:** *Let F be a field of characteristic zero. Let A be a proper subset of F which is a semifield under the operations of F. Then we say F to be a Smarandache definite special field. (S-definite special field).*

We illustrate this by the following example.

*Example 3.1.17:* Let Q be the field of characteristic zero. Take $A = Z^+ \cup \{0\}$, A is a semifield and $A \subset Q$ so Q is a S-definite special field.

Now we define Smarandache definite special subfield.

**DEFINITION 3.1.8:** *Let F be any field. Let K be a subfield of F, we call K to be Smarandache definite special subfield (S-definite special subfield) if K itself is a S-definite special field.*

Thus we see if F is a field having a subfield which is a S-definite special subfield then F itself is a S-definite special field.

**THEOREM 3.1.21:** *Let F be any field. K a proper subset of F which is a subfield and K is a S-definite special subfield of F; then F itself is a S-definite special field.*

*Proof:* Given F is a field and K a proper subset of F which is a subfield of F which is a subfield and K is also given to be a S-definite special subfield of F. To show F is itself a S-definite special field. Since K is given to be a S-definite special subfield; K contains a proper subset $T \neq \phi$ ($T \subseteq K$) such that T is a semifield. Now $K \subseteq F$ and $T \subseteq K \subset F$ i.e., $T \subseteq F$ and T is a semifield so F is a S-definite special field. Hence the claim.

We illustrate this situation by the following example.

*Example 3.1.18:* Let R be the field of reals. R contains a proper subset Q and Q is the field of rationals. We see Q is a S-definite special field as Q contains the subset $Z^+ \cup \{0\} = V$, where V is a semifield. Now $V \subseteq Q$ so Q is a S-definite special subfield of R. Now $Q \subseteq R$ and $V \subseteq Q \subseteq R$ and V is a subset of R which is a semifield. Hence R is a S-definite special field.

Now however the converse of the above statement in general is not true which is given by the following theorem.

**THEOREM 3.1.22:** *Let F be any field which is a S-special definite field. F in general need not contain a S-definite special subfield.*

*Proof:* We prove this only by a counter example. Take Q to be the field of rationals. Clearly Q is a S-definite special field for Q contains $V = Z^+ \cup \{0\}$ which is a semifield. Now Q being the field of rationals does not contain any proper subfields so Q does not contain any S-definite special subfields. Hence the claim.

We call all fields F which do not contain S-definite special subfields and which is a S-definite special field as S-definite special prime field. Q is a S-definite special prime field.
R is not a S-definite special prime field.

For more above semifields please refer [102].

### 3.2 Smarandache Special Vector Spaces

Now having defined the notion of S-definite special fields we now proceed on to define the new notion of Smarandache definite special vector spaces and give a few important properties about them.

**DEFINITION 3.2.1:** *Let F be a field which is a S-definite special field and V be a vector space over F. We say V is Smarandache definite special vector space (S-definite special vector space) over F if V has a proper subset W such that W is a semivector space over at least one semifield contained in F.*

We illustrate this by some examples.

*Example 3.2.1:* Let Q[x] be a vector space over Q. We see Q is a S-definite special field. Let $Z° = Z^+ \cup \{0\}$ be the semifield in Q. Let Z°[x] be the proper subset of Q[x]. Clearly Z°[x] is a semivector space over Z°, so Q[x] is a S-definite special vector space over Q. In fact we can take $Q° = Q^+ \cup \{0\}$ to be a proper subset and $Q^o[x]$ is also a semivector space over $Z^o$, so $Q[x]$ is a S-definite special vector space over Q.

***Example 3.2.2:*** Let $M_{n \times n} = \{M = (m_{ij}) \mid m_{ij} \in Q\}$ be the vector space over Q. $M_{n \times n}$ is a S-definite special vector space over Q, because Q is a S-definite special field with $Z^\circ = Z^+ \cup \{0\}$ and if we choose $P_{n \times n} = \{P = (p_{ij}) \mid p_{ij} \in Z^\circ\}$ then $P_{n \times n}$ is a semivector space over $Z^\circ$. Hence the claim.

Now we proceed on to define the new notion of Smarandache definite special basis of a S-special definite vector space.

**DEFINITION 3.2.2:** *Let V be a vector space over F. Suppose V be a S-definite special vector space over F and if $B = \{x_1, ..., x_n\}$ forms a basis of a semivector space in V as well as the vector space V over F then we call a B to be S-definite special basis of V.*

We illustrate this situation by the following example.

***Example 3.2.3:*** Let $Q \times Q \times Q = V$ be a vector space over Q. V is a S-definite special vector space. For take $T = Z^\circ \times Z^\circ \times Z^\circ$; T is a semivector space over the semifield $Z^\circ$. $\{(1\ 0\ 0), (0\ 1\ 0), (0\ 0\ 1)\}$ is a basis of T as well as a basis of V, so T is a S-definite special basis of V. Take $B = \{(0\ 3\ 0), (0\ 0\ 1), (4\ 0\ 0)\}$ to be a basis of V. Clearly B is not a S-definite special basis of V for B is not a basis of $T = Z^\circ \times Z^\circ \times Z^\circ$.

We give yet another example.

***Example 3.2.4:*** Let $Q[x]$ be a vector space over Q. $Z^\circ[x]$ is a semivector space over $Z^\circ$, the semifield. $\{1, x, x^2, ...\}$ is a basis of $Z^\circ[x]$ over $Z^\circ$ as well as $Q[x]$ over Q i.e., $\{1, x, x^2, ..., x^n, ...\}$ is a S-definite special basis of $Q[x]$. Take $\{1/3, x, x^2, x^2, ..., x^n, ...\} = B$. B is a basis of $Q[x]$ but is not a basis of $Z^\circ[x]$ so B is not a S-definite special basis of $Q[x]$.

We now define the S-definite special dimension of a S-definite special vector space.

**DEFINITION 3.2.3:** *Let V be a vector space over F. S $\subseteq$ V be a semivector space over T $\subset$ F. If dimension of V is the same as dimension of at least one S $\subset$ V (say) n, then we say the Smarandache definite special dimension (S- definite special dimension) of the vector space is n where n is the number of elements in the basis of the vector space V and n is also the number of elements in at least one of the semivector spaces S $\subseteq$ V over the semifield T $\subset$ F. If one does not get such n we say the S-definite special dimension of the S-definite special vector space is not defined.*

We illustrate this by the following example.

*Example 3.2.5:* Let V = R × R × R × R be a vector space over the field R. We see V as a vector space over R is of dimension 4. W = R° × R° × R° × R° (where R° = R $\cup$ {0}) is a semivector space over Q° = Q$^+$ $\cup$ {0} of infinite dimension. W as a semivector space over R° = R$^+$ $\cup$ {0} is of dimension four and {(1, 0, 0, 0), (0, 1, 0, 0), (0, 0, 1, 0), (0, 0, 0, 1)} = B is a basis of W over R°. Now B is also a basis of V. Hence the S-definite special basis of V over R is four.

Now W as a semivector space over Q° is infinite dimensional. Thus we see in the first instance one can see that the number of elements in the basis of the semivector space is not always equal to the number of elements in the basis of the given vector space V over R. Thus the S-definite special dimension is given by the existence of atleast one semivector space in V, having the same number of base elements. Hence the S-definite special dimension in this case is however four. We see we can find several semivector spaces over appropriate semifields contained in the given field.

We illustrate this by the following example.

*Example 3.2.6:* Consider the vector space $M_{3\times 2}$ = {$(m_{ij})$ = M | $m_{ij}$ $\in$ Q, 1 $\leq$ i $\leq$ 3, 1 $\leq$ j $\leq$ 2} over Q. Clearly $M_{3\times 2}$ is a S-definite special vector space as $M_{3\times 2}$ contains $P_{3\times 2}$ = {p = $(p_{ij})$ | $p_{ij}$ $\in$ Z° =

$Z^+ \cup \{0\}$; $1 \le i \le 3$ and $1 \le j \le 2\}$ to be a proper subset which is a semivector space over $Z^\circ$ clearly a basis of $P_{3\times 2}$ is six given by

$$B = \left\{ \begin{pmatrix} 1 & 0 \\ 0 & 0 \\ 0 & 0 \end{pmatrix}, \begin{pmatrix} 0 & 1 \\ 0 & 0 \\ 0 & 0 \end{pmatrix}, \begin{pmatrix} 0 & 0 \\ 1 & 0 \\ 0 & 0 \end{pmatrix}, \begin{pmatrix} 0 & 0 \\ 0 & 1 \\ 0 & 0 \end{pmatrix}, \begin{pmatrix} 0 & 0 \\ 0 & 0 \\ 1 & 0 \end{pmatrix}, \begin{pmatrix} 0 & 0 \\ 0 & 0 \\ 0 & 1 \end{pmatrix} \right\}.$$

B is also a basis of the vector space $M_{3\times 2}$ over Q. Thus

1. $M_{3\times 2}$ is a S-definite special vector space.
2. This S-definite special vector space has B to be a S-definite special basis.
3. The S-definite special dimension of $M_{3\times 2}$ over Q is six.

Now consider the proper subset $S_{3\times 2} = \{(s_{ij}) = S \mid 1 \le i \le 3, 1 \le j \le 2, s_{ij} \in Q^0 = Q^+ \cup \{0\}\}$ in $M_{3\times 2}$. $S_{3\times 2}$ is a semivector space over $Z^\circ = Z^+ \cup \{0\}$. Also $S_{3\times 3}$ is a semivector space over $Q^\circ = Q^+ \cup \{0\}$. Now the number of basis of $S_{3\times 2}$ as a semivector space over $Z^\circ = Z^+ \cup \{0\}$ is infinite. But $S_{3\times 2}$ as a semivector space over $Q^\circ = Q^+ \cup \{0\}$ is finite and it is six. It is interesting to note that $P_{3\times 2}$ is not a semivector space over $Q^\circ = Q^+ \cup \{0\}$ as it is evident from the fact any $p/q \in P_{3\times 2}$ ($q \ne 0$) $p \ne q$. Hence the claim.

We give yet another example in this direction.

***Example 3.2.7:*** Let $V = R \times R$ a vector space over Q. V is S-definite special vector space as $W = Z^\circ \times Z^\circ \subseteq V$ is a semivector space over $Z^\circ \subseteq Q$. But the dimension of V is infinite over Q and dimension of W over $Z^\circ$ is two given by the basis $B = \{(0, 1) (1, 0)\}$. Now what is the S-definite dimension of V as a vector space. To this end we must try to find a semivector space of infinite dimension.
    Now take $P = R^\circ \times R^\circ$ where $R^\circ = R^+ \cup \{0\}$, P is a semivector space over $Q^\circ = Q^+ \cup \{0\} \subseteq Q$. This P is infinite dimensional as it has infinite basis. Thus any basis of P will be a basis of V also. Now consider $T = Q^\circ \times Q^\circ$, $T \subseteq V$ and T is a semivector space over semifield $Q^\circ$ and $B^1 = \{(1,0), (0,1)\}$ is a

basis of T over Q°. Also we see certainly $B^1$ is not a basis of the vector space V over Q. Now the same semivector space T is defined over the semifield Z°. Now the basis of T relative to Z° is infinite but that infinite is different from the infinite basis of V as a vector space over Q.

Infact this infinite set which acts as a basis for T over Z° cannot act as a basis of V over the field Q. Thus it is not a S-definite special basis. Now consider R° × R° = M; M is a semivector space over Z° and the number of elements in the basis is also infinite but this is yet another infinity different from the infinite basis of semivector space T over Z°. Thus we have seen in case of S-definite special vector spaces the semivector spaces in them are many and with varying dimension.

We suggest some problems about them.

Next we proceed on to define the notion of Smarandache definite special linear transformation of Smarandache special definite vector spaces.

**DEFINITION 3.2.4:** *Let V and W be vector spaces over the same field F we say T from V to W to be a Smarandache definite special linear transformation (S-definite special linear transformation) if we have T to be a semivector space linear transformation of a semivector space S from V to $S_1$ from W over the same semifield $K \subset F$ i.e., T: $S \to S_1$ such that T (au + v) = aT(u) + T(v) for all u, v $\in$ S and a $\in$ K.*

*Note:* It is important to note that we can have several linear transformations with different domain and range spaces for the same S-definite special vector spaces V and W. Further it is made very clear that all of the domain space need not be transformed under T, infact T when we call it as a S-definite special linear transformation we only make it work on the semivector space S contained in V and not on the whole of V.

We illustrate this by the following example.

*Example 3.2.8:* Let T**:** V $\to$ W be a S-definite special linear transformation where V = Q × Q × Q × Q × Q defined over Q

and W = R × R × Q × Q defined over Q. Here T is from the semivector space S = Z° × Z° × Z° × Z° × Z° over Z° and $S_1$ = Q° × Q° × Z° × Q° over Z°. T(x, y, z, u, v) = (x + y, u, y, z + u + v); it is easily verified T is a S-definite special linear transformation of the two S-definite special vector spaces V and W both defined over Q. Suppose we take S = Z° × Z° × Q° × Q° × Q° over Z° and $S_1$ = Q° × Q° × Q° × Q° over Q° then clearly we cannot define a S-definite special linear transformation form V to W restricting it to S and $S_1$ respectively (S ⊂ V and $S_1$ ⊂ W).

Several interesting properties regarding the S-definite special linear transformation can be obtained. In this direction the authors have suggested several problems in the last chapter of this book.

Now we proceed on to define the notion of Samarandache definite special subvector space of a vector space.

**DEFINITION 3.2.5:** *Let V be any vector space over the field F. Suppose W be a subspace of V then we say W is a Smarandache definite special subvector space (S-definite special subvector space) of V if W has a proper subset P such that P is a semivector space over a semifield K where K is a proper subset of F.*

We illustrate the definition by the following example.

*Example 3.2.9:* Let V = Q × Q × Q be a vector space over Q. Let W = Q × Q × {0}, be a proper subset of V, we prove W is a Smarandache special definite subspace of V. Take the proper subset T = Z° × Z° × {0} of W. Clearly T is a semivector space over Z° where Z° ⊆ Q. This W is a S-special definite vector subspace of V.

In view of the definition of S-definite special subvector space we prove the following theorem.

**THEOREM 3.2.1:** *Let V be a vector space over a field F. Suppose W ⊂ V is a proper subset of V which is a S-definite special subspace of V then V is itself a S-definite special vector space.*

*Proof:* Given V is a vector space over F and W is a proper subset of V which is a S-definite special subspace of V. Now W is a S-definite special vector space over F so W has a proper subset T such that T is a semivector space over the semifield K where K is a proper subset of F. Now T ⊂ W but W ⊂ V so T ⊂ W but W ⊂ V so T ⊂ W ⊂ V i.e., T is a proper subset of V i.e. T ⊆ V and T is a semivector space over K; K a semifield contained in F so V is itself a S-definite special vector space. Hence the claim.

**THEOREM 3.2.2:** *Let V be a vector space over the field F. V is a S-special definite vector space over F then V is a S-special definite group and F is a S-special definite field.*

*Proof:* Suppose V is given to be a S-definite special vector space over F to show V is a S-special definite group and F is a S-definite special field. Now V is a S-definite special vector space implies V has a proper subset S such that S is a semivector space over the semifield K, K a proper subset of F which is a semifield. We know (V, +) is an ableian group by the very definition of a vector space.

Now S ⊂ V and (S, +) is a commutative semigroup. Hence (V, +) is a S-definite special group. Further we know since V is a S-definite special vector space over F. We have K ⊂ F and K is a semifield. Thus when a field F contains a proper subset K, which is a semifield then we know F is a S-definite special field. Hence one way is true.

Now on the other hand we are given (V, +) is a S-definite special group so we can find in V a proper subset W such that (W, +) is a semigroup. Now we are also given F is a S-definite special field so F contains a proper subset K such that K is a semifield. We know FV ⊂ F i.e., for all v ∈ V and a ∈ F, av ∈ F so for all w ∈ W and a ∈ K, aw ∈ V but we need aw to be in W but this may happen or may not happen in general. Thus the converse of the theorem may not be true in general.

This we illustrate by the following example.

***Example 3.2.10:*** Let $V = Q \times Q \times Q$ be a S-definite special group with respect to addition. Q be a S-definite special field. Suppose $W = Z° \times Z° \times Z°$ be a semigroup in V and $Q°$ be the semifield in Q we see for any $a \in Q°$ and $w \in W$ $aw \notin W$; for take

$$a = \frac{3}{7}, w = (1\ 1\ 0);\ aw = (\frac{3}{7}, \frac{3}{7}, 0) \notin W = Z° \times Z° \times Z°.$$

Hence the claim.
   This we see the converse of the theorem 3.2.2 is not in general true.

**DEFINITION 3.2.6:** *Let V be a S-definite special vector space over F. We call a function ( | ) on a semivector space $W \subseteq V$ which assigns to each ordered pair of vectors α, β $\in$ W a scalar (α | β) in K, K a semifield in F in such a way that for all α, β, γ in W and for all scalars c in K.*

   *(a) (α + β | γ) = (α|γ) + (β|γ).*
   *(b) (cα | β) = c(α | β).*
   *(c) (β | α) = $\overline{(α | β)}$ if elements are from complex semivector space.*
   *(d) (α | α) > 0 if α ≠ 0.*

*Then we call the function ( | ) to be a Smarandache definite special inner product on V (S-definite special inner product on V).*

*Note:* It is important to note that the S-definite special inner product on V need not be defined on whole of V if it is defined on a semivector space W in V it is sufficient. So inner product also need note be defined on the whole of V. As usual if V = Q

× Q × Q then an inner product on V, for x = ($x_1$, $x_2$, $x_3$) and y = ($y_1$, $y_2$, $y_3$) ∈ W = Z° × Z° × Z° defined by (x | y) = $\sum x_i y_i$ will be called as S-definite special standard inner product.

As in case of other vector spaces one in case of S-definite special inner product defined on V cannot claim |(α|β)| < ||α|| ||β||.

Further writing

$$\beta = \frac{(\beta/\alpha)}{\|\alpha\|^2} \alpha$$

may not be possible for all vector space V for we see in case the semivector space of V is of the form W ⊂ V, W = Z° × Z° × … × Z° then every element in the n-tuple can only be positive integers hence the claim; i.e., we wish to bring to the notice of the readers that for every inner product defined on the vector space all its properties cannot be restricted to the subset which is the semivector space of the vector space. Thus one cannot blindly say all the properties satisfied by the inner product on V is true in case of the S-definite special inner product on V. With this in mind the reader can as a researcher develop properties about S-definite special inner products on V.

We cannot define the notion of orthogonal vectors in general using the notion of S-special definite standard inner product or we can show that there exist S-definite special inner products which have no element other than zero to be orthogonal. So while characterizing S-definite special vector spaces we see the Gram Schmidt process cannot be applied for all S-definite special inner products on S-definite special vector spaces.

## 3.3 Smarandache Definite Special Linear Algebra

In this section we proceed on to define the new notion of Smarandache definite special linear algebra.

**DEFINITION 3.3.1:** *Let V be a semivector space over a semifield F. We say V is a semilinear algebra if for each pair of vectors v, u in V vu and uv are in V; called product of u and v in such a way that*

a. *Multiplication is associative; u(vt) = (uv)t for all u, v, t ∈ V.*
b. *Multiplication is distributive with respect to addition;*
    *u(v + t ) = uv + ut and*
    *(u + v)t = ut + vt for all u, v, t ∈ V.*
c. *For each scalar c in F c(uv ) = (cu )v = u( cv ) and u, v ∈ V.*
d. *If there is an element 1 in V such that u.1 = 1.u = u for each u ∈ V;*

*we call V to be the semilinear algebra with identity over F. If in V for every u, v ∈ V, uv = vu, we call V to be a commutative semilinear algebra.*

We first illustrate this by some examples before we proceed to define the notion of S-definite special linear algebra.

**Example 3.3.1:** $Z^o \times Z^o \times Z^o = V$ be a semivector space over the semifield $Z^o$. Clearly V is a semilinear algebra by defining for x = $(x_1, x_2, x_3)$ and y = $(y_1, y_2, y_3)$ in V we have x.y = $(x_1y_1, x_2y_2, x_3y_3)$ ∈ V. Hence the claim.

**Example 3.3.2:** Let $M_{n \times n}$ = {M = $(m_{ij})$ | $m_{ij} \in Z^o$} be the collection of all n × n matrices with entries from $Z^o$, that is for A, B ∈ $M_{n \times n}$, A · B ∈ $M_{n \times n}$; where "·" is the matrix multiplication. $M_{n \times n}$ is a semilinear algebra over $Z^o$.

Now having seen the definition and examples of semilinear algebra, we proceed on to define the notion of Smarandache definite special linear algebra.

**DEFINITION 3.3.2:** *Let V be a linear algebra over a field F. We say V is a Smarandache definite special linear algebra (S-definite special linear algebra) if V contains a nonempty subset W in V such that W is a semilinear algebra over a semifield K contained in F.*

***Example 3.3.3:*** Let $V = Q \times Q \times Q \times Q$ be a linear algebra over Q. Take $W = Z^o \times Z^o \times Z^o \times Z^o$ over $Z^o \subseteq Q$ where $Z^o = Z^+ \cup \{0\}$; W is a S- definite special linear algebra over V, as W is a semilinear algebra over $Z^o \subseteq Q$.

We see all S-definite special vector spaces need not be S-definite special linear algebras.
  We can only illustrate this situation by examples.

***Example 3.3.4:*** Let $M_{3\times 2} = \{M = (m_{ij})|\ m_{ij} \in Q \}$ be the collection of all $M_{3\times 2}$ matrices with entries from Q. $M_{3\times 2}$ is only a S- definite special vector space over Q and is not a S-definite special linear algebra over Q. For $M_{3\times 2}$ does not contain any semilinear algebra over a semifield.
  Thus we can say that as in case of vector spaces, all linear algebras are vector spaces but vector spaces in general are not linear algebras; like wise all S-definite special linear algebras are S-definite special vector spaces but S-definite special vector spaces in general need not be S-definite special linear algebras.
  We now illustrate by examples some linear algebras which are not S-definite special linear algebras.

***Example 3.3.5:*** Let $V = Z_2 \times Z_2 \times Z_2$ be a linear algebra over $Z_2$. Clearly V is not a S-definite special linear algebra over $Z_2$.

***Example 3.3.6:*** Let $M^2_{3\times 3} = \{M = (m_{ij}) \mid m_{ij} \in Z_7\}$ be the collection of all $3 \times 3$ matrices with entries from $Z_7$. $M^2_{3\times 3}$ is a linear algebra over $Z_7$. We see $M_{3\times 3}$ is not a S-definite special linear algebra as $Z_7$ does not contain any proper subset which is a semifield.

Now we will define the notion of Smarandache definite special linear algebra.

**DEFINITION 3.3.3:** *Let V be a linear algebra over a field F. If W is a proper subset of V and W is a S- definite special linear algebra then we call W to be a S- definite special linear subalgebra of V over F.*

The following claim given in the form of the theorem which can be easily proved.

**THEOREM 3.3.1:** *Let V be a linear algebra over the field F. If V has a proper subset W such that W is a S-definite special sub linear algebra over F then V is itself is a S-definite special linear algebra over F.*

We illustrate this situation by the following example.

*Example 3.3.7:* Let $V = Q \times Q \times Q \times Q \times Q$ be a linear algebra over Q. Let $W = Q \times Q \times Q \times \{0\} \times \{0\}$ be a proper subset of V. Clearly W is a S- definite special linear algebra over Q. For take $T = Z^o \times Z^o \times Z^o \times \{0\} \times \{0\}$, a proper subset of W, we see T is a semilinear algebra over $Z^o$. Hence W is a S- definite special linear algebra which is nothing but a S-definite special linear subalgebra of V. Now we claim V is itself a S-definite special linear algebra, for $T = Z^o \times Z^o \times Z^o \times \{0\} \times \{0\} \subseteq W \subseteq V$; so T contained in V and T is a semilinear algebra over $Z^o$, hence V is a S-definite special linear algebra over Q.

Now we have seen some of the properties of S-special linear algebra. We show the existences of a few class of S-definite special linear algebras; which is stated by the following theorems.

**THEOREM.3.3.2:** *Let F[x] be the linear algebra of polynomials over the field F of characteristic zero then, F[x] is a S-special definite linear algebra over F.*

*Proof:* Given F[x] is a linear algebra over the field F of characteristic zero. Since F is a field of characteristic zero we see either Q the field of rationals is a subfield of F or F is isomorphic with Q. Now in both cases we see F contains nontrivial semifields like $Q^o = Q^+ \cup \{0\}$ and $Z^o = Z^+ \cup \{0\}$. Consider the collection of all polynomials with coefficients from $Q^o$ denoted by $Q^o[x]$. $Q^o[x]$ is clearly a semilinear algebra over the semifield Q also $Q^o[x]$ is a semilinear algebra over the semifield $Z^o$. As $Q^o[x] \subseteq F[x]$, we see $F[x]$ is a S-special definite linear algebra over F. Hence the theorem.

It is pertinent to mention here that $Z_p[x]$ where $Z_p$ is a prime field of characteristic p (p a prime) and $Z_p[x]$ the collection of all polynomials with coefficients in $Z_p$ is not a S- special definite linear algebra over $Z_p$.

In view of this we propose some problems in the last chapter of this book.

*Remark:* Suppose we restrain the degree of polynomials to be less than or equal to say n, that is $x^{n+2} = x$ or $x^{n+1} = 1$, then we by replacing all polynomials p(x) of degree m in x if $m \leq n$ then we take p(x) as the polynomial. If $m > n$ then we replace those powers of x which are greater than n by m – n.

(If the reader has some confusion over this he/she can make use of this illustrative example.) Suppose
$$F[x] = \{\sum_{i=0}^{5} a_i x^i \mid x^6 = x^0 \text{ and } a_i \in F\}.$$
Then F[x] is clearly a vector space over F and dimension of F[x] over F is 6 and the standard basis is given by $\{1, x, x^2, x^3, x^4, x^5\}$. Now we can treat or make F[x] as a linear algebra over F by declaring multiplication in this way; replace all polynomials after usual polynomial multiplication if it contains powers of x greater than 5 then put $x^6 = x^0 = 1$.
For instance $x^{13} = x$, $x^7 = x$ and $x^9 = x^3$. Suppose $p(x) = x^5 + 3x^3 - 4x + 1$ and $q(x) = 3x^5 + 2x^2 - 7$ then $p(x)q(x) = (x^5 + 3x^3 - 4x + 1)(3x^5 + 2x^2 - 7) = 3x^{10} + 9x^8 - 12x^6 + 3x^5 + 2x^7 + 6x^5 - 8x^3 + 2x^2 - 7x^5 - 21x^3 + 28x - 7$ by putting $x^6 = x^0 = 1$, we get $p(x)q(x) = 2x^5 + 3x^4 - 29x^3 + 11x^2 + 30x - 19$.

Clearly $p(x)q(x) \in F[x]$. Thus F[x] is a linear algebra of finite dimension, infact F[x] is a finite dimensional S- definite special linear algebra. We have given a class of S-definite special linear algebra of both finite and an infinite dimension having standard basis $\{1, x, x^2, x^3, …, x^n \mid n < \infty\}$ and $\{1, x, x^2, …, x^n, …\}$ respectively.

We give yet another class of S- definite special linear algebras by this theorem.

**THEOREM 3.3.3:** *Let $M_{n \times n} = \{M = (m_{ij}) \mid m_{ij} \in F; 1 \le i, j \le n\}$ be the set of all $n \times n$ matrices with entries from the field $F$ of characteristic zero. $M_{n \times n}$ is a S- definite special linear algebra.*

*Proof:* Given F is a field of characteristic zero so the field of rationals is either a subfield of F or is isomorphic with Q. It is further given $M_{n \times n}$ is a linear algebra over F. Now to show $M_{n \times n}$ is a S-definite special linear algebra over F. We see F is a S-definite special field as F contains $Z^o = Z^+ \cup \{0\}$ and $Q^o = Q^+ \cup \{0\}$ as semifields. Now $P_{n \times n} = \{P = (p_{ij}) \mid p_{ij} \in Q^o = Q^+ \cup \{0\}, 1 \le i, j \le n\}$. Clearly $P_{n \times n} \subseteq M_{n \times n}$ and $P_{n \times n}$ is semilinear algebra over the semifield $Z^o$. Also $P_{n \times n}$ is a semilinear algebra over the semifield $Q^o = Q^+ \cup \{0\}$. This $M_{n \times n}$ is a S- definite special linear algebra. Infact $M_{n \times n}$ is a S- definite special linear algebra over F.

Thus we have given yet another class of S- definite special linear algebra.

*Example 3.3.8:* Let $M_{2 \times 2} = \{M = (m_{ij}) \mid m_{ij} \in Q; 1 \le i, j \le 2\}$ be the set of all $2 \times 2$ matrices with entries from Q. $M_{2 \times 2}$ is a S-definite special linear algebra as $P_{2 \times 2} = \{(p_{ij}) = P \mid p_{ij} \in Z^o = Z^+ \cup \{0\}; 1 \le i, j \le 2\}$ is a semilinear algebra over the semifield $Z^o = Z^+ \cup \{0\}$. Thus $M_{2 \times 2}$ is a S- definite special linear algebra over Q. Take

$$B = \left\{ \begin{pmatrix} 1 & 0 \\ 0 & 0 \end{pmatrix}, \begin{pmatrix} 0 & 1 \\ 0 & 0 \end{pmatrix}, \begin{pmatrix} 0 & 0 \\ 1 & 0 \end{pmatrix}, \begin{pmatrix} 0 & 0 \\ 0 & 1 \end{pmatrix} \right\};$$

B is a S-definite special basis of the S-definite special linear algebra as B is a basis of $P_{2 \times 2}$ as a semilinear algebra over the semifield $Z^o$ as well as B is a basis of $M_{2 \times 2}$, the linear algebra over Q. Hence the claim.

Another class of S-definite special linear algebra is given by the following theorem.

**THEOREM 3.3.4:** *Let* $V = \underbrace{Q \times Q \times \ldots \times Q}_{n-times}$ *be a linear algebra over Q. V is a S-definite special linear algebra of S-definite dimension n.*

*Proof:* Given $V = Q \times Q \times Q \times \ldots \times Q$; n-times is a linear algebra over Q. Take $S = Q^o \times Q^o \times \ldots \times Q^o$; n-times contained in V where $Q^o = Q^+ \cup \{0\}$, clearly S is a semilinear algebra over $Q^o$, the semifield in Q. Hence V is a S-definite special linear algebra over Q. Now $B = \{(1, 0, 0, \ldots, 0), (0, 1, 0, \ldots, 0), \ldots, (0, 0, 0, \ldots, 1)\}$ is a basis of both V and S, hence B is a S-definite special basis of the S-definite special linear algebra over Q.

We give the following example to the interested reader.

*Example 3.3.9:* Let $V = Q \times Q \times Q \times Q$ be a linear algebra over Q. V is a S-definite linear algebra over Q as V contains a proper subset $T = Z^o \times Z^o \times Z^o \times Z^o \subseteq V$ and $Z^o \subseteq Q$ so T is a semilinear algebra over $Z^o$, the semivector space contained in Q. Hence T is a S-definite special linear algebra over Q. Clearly $B = \{(1, 0, 0, 0), (0, 1, 0, 0), (0, 0, 1, 0), (0, 0, 0, 1)\}$ is a S-definite special basis of the S-definite special linear algebra over Q.

Now we proceed on to define the notion of Smarandache definite special transformation of S-definite special linear algebra.

**DEFINITION 3.3.4:** *Let V and W be any two S-definite special linear algebras defined over the fields F and K respectively both of characteristic zero. We call a linear transformation from V to W to be a Smarandache definite special linear transformation (S-definite special linear transformation) if T is a linear transformation from a semilinear algebra M to N, where M is a semilinear algebra in V and N is a semilinear algebra in W; but* *both M and N are semilinear algebra defined on the same semifield P contained in both F and K, that is; T:V → W is such that T: M (⊆ V)→ N(⊆ W), that is T need not be defined on the whole of V but only defined on the semilinear algebra M contained in V to the semilinear algebra N contained in W, but both M and N are semilinear algebras over the same semifield P.*

The following are the difference between linear transformation of the linear algebras and the S-definite special linear algebras.

1. T in case of linear transformation defined on whole of V to W and T in case of S-definite special linear transformation defined only on the semivector space contained in V need not be defined on whole of V.

2. In case of linear transformation we need both V and W to be defined as linear algebras over the same field, but incase of S-definite special linear transformation V and W the S-definite special linear algebras need not be defined over the same field but only the semilinear algebra M and N contained in the linear algebras V and W respectively be defined on the same semifield.

3. The linear algebra transformation from V to W in general even if restricted from M (M ⊆ V) to N (N ⊆ W) need not be well defined. If in case the linear transformation be defined we will give a new name for such transformations.

We first illustrate these situations by the following examples.

*Example 3.3.10:* Let $V = Q \times Q \times Q \times Q$ be a S- definite special linear algebra defined over the field Q and let $W = R \times R$ be a S- special definite linear algebra defined over R. We call a map T: V → W to be a S-definite special linear algebra transformation, that is; T: M ⊆ V → N ⊆ W where T(x, y, z, w) = (x + y, z + w) for all x, y, z, w ∈ $Q^o \times Q^o \times Q^o \times Q^o$ = M, a semilinear algebra over the semifield $Q^o$ where $Q^o = Q^+ \cup$

{0} and $M \subseteq V$; here $N = Q^o \times Q^o \subseteq W$ a semilinear algebra over the semifield $Q^o$.

We see any linear transformation of linear algebras in general need not be S-definite special linear algebra linear transformations even when restricted to semilinear algebras. This would be made explicit by the following example.

*Example 3.3.11:* Let $V = Q \times Q \times Q \times Q$ be a S-definite special linear algebra over Q. Let $W = R \times R \times R \times R \times R$ be a S-definite special linear algebra over R. Let $M = Z^o \times Z^o \times Z^o \times Z^o$ defined over $Z^o$ be a semilinear algebra contained in V and $N = Q^o \times Q^o \times Q^o \times Q^o \times Q^o$ be a S-definite special linear algebra defined over $Z^o$, N a semilinear algebra contained in W. Let $T: V \to W$ be a S-definite linear transformation of S-definite special linear algebras from M to N defined by $T(x, y, z, w) = (x + y, y, z, w + y, z)$ for all $(x, y, z, w) \in Z^o \times Z^o \times Z^o \times Z^o \subseteq V$ and $(x + y, y, z, w + y, z) \in N = Q^o \times Q^o \times Q^o \times Q^o \times Q^o \subseteq W$. Clearly T is a S-definite special linear algebra transformation of S-definite special linear algebras.

*Example 3.3.12:* Let V and W be S-definite special linear algebras defined over the S-definite special field Q, where $V = Q \times Q \times Q \times Q \times Q$ and $W = Q \times Q \times Q$ defined over Q. Let $T: V \to W$ be a linear transformation given by $T(x, y, z, w, u) = (x - y, - y, - z + x)$ for $(x, y, z, w, u) \in V$ and $(x - y, - y, -z -x) \in W$.

We see T is a linear algebra linear transformation. Clearly T is not a S-definite linear transformation of S-definite special linear algebras.

Yet we give another example before we give our conclusions.

*Example 3.3.13:* Let V and W be any two S-definite special linear algebras where $V = R \times R \times R$ defined over Q and $W = Q \times Q \times Q \times Q$ defined over Q. Consider $P = Z^o \times Z^o \times Z^o \subseteq V$ be a semilinear algebra over $Z^o$ and $S = Z^o \times Z^o \times Z^o \times Z^o \subseteq W$ be a semilinear algebra over $Z^o$.

Define S- definite linear transformation T: V $\to$ W by T(x, y, z) = (x + y, z, y + z, x). T is a S-definite special linear algebra transformation.

Can T be extended to V so that T is a linear algebra transformation from V to W?

Now we give an application which becomes an important one in case of defining Smarandache definite Leontief economic models. [104]

We see suppose for any Smarandache definite vector space of n×n matrices over the field F of characteristic zero.

Suppose we have S- definite vector space V of n×n matrices with entries from the real field then if we take any map from the S-definite vector space V to the semivector space over $R^o$ say P = {A = $(a_{ij})$ | $a_{ij} \in R^o$} $\subseteq$ V.

Then we see any Leontief model with n entries belong to P, so these models can be easily generated by this particular S-definite vector space over reals.

To this end we define a new type of linear transformation on S- definite vector spaces.

**DEFINITION 3.3.5:** *Let V be any S-definite vector space over the field K of characteristic zero. Let W $\subseteq$ V be any non empty proper subset of V such that W is a semivector space over some semifield F $\subseteq$ K. A linear transformation $T_s(x + y) = T_s(x) + T_s(y)$ and $T_s(ax) = a\, T_s(x)$ for all x, y $\in$ V and for all a $\in$ K will be known as Smarandache definite special converging linear transformation (S-definite special converging linear transformation ) on V.*

We first illustrate it by an example.

***Example 3.3.14:*** Let V = Q × Q × Q × Q be a S- definite special vector space over Q. Let W = $Z^o \times Z^o \times Z^o \times Z^o$ be a semivector space over $Z^o$. Let T: V $\to$ W be defined by T(x, y, z, w) = T (|x|, |y| + |z|, |z| = |w|, |w|). Then it is easily verified that T is a S-special definite converging linear transformation on V.

All identity linear transformation of V to V are S-definite special converging linear transformation when restricted to any of the semivector spaces contained in V. One may wonder can we ever have the concept of Smarandache special definite diverging linear transformation of a vector space V? The answer is yes and we define a Smarandache special definite diverging linear transformation.

**DEFINITION 3.3.6:** *Let V be a S-definite special vector space over a field F of characteristic zero. Let T be a map from V to V that is; T is a map from a semivector space $W \subseteq V$ to V such that $T(ax + y) = aT(x) = T(y)$ for all $x, y \in W$ and $a \in F$. Such a map T is defined as the Smarandache special definite diverging linear transformation (S-special definite diverging linear transformation) of V.*

We first illustrate this situation by the following example.

*Example 3.3.15:* Let $V = R \times R \times R \times R$ be a S-definite special vector space over R. Let $T: W \subseteq V \to V$ be a S- definite special diverging linear transformation of V where $W = R^o \times R^o \times R^o \times R^o$ is a semivector space over $R^o$ defined by $T(x, y, z, w) = (x - y, y - z, z - w, w - x)$ for all $(x, y, z, w) \in W$.

In general a S-definite special divergence linear transformation need not be a S- definite special convergence linear transformation.

Next we proceed on to define the notion of Smarandache definite special near rings and develop some important properties about these S-definite special near rings.

**DEFINITION 3.3.7:** *Let N be a near ring, we say N is a Smarandache definite special near ring (S-definite special near ring) if N has a proper subset S which is a seminear ring under the operations of N and $\phi \neq S \subseteq N$.*

We see all near rings are seminear rings but near rings in general are not seminear rings. We see the class of seminear rings strictly contains the class of near rings; that is, seminear rings are the generalization of near rings, i.e., in Smarandache definite special algebraic structures we demand a weaker structure to be contained in the stronger structure where as the Smarandache algebraic structure are those algebraic structures which contain in them a subset which is a stronger structure imbedded in the weaker structure.

First we illustrate this by the following example.

***Example 3.3.16:*** Let $\{Z, +, \odot\}$ be a near ring where + is a group under + and $\odot$ be a operation on Z defined by $a \odot b = a$ for all $a, b \in Z$. $\{Z^o, +, \odot\}$ is a near ring. We see $\{Z^o, +, \odot\}$ is a proper subset of Z and $Z^o$ is a seminear ring. Hence $\{Z, +, \odot\}$ is a S-definite special near ring.

Now we give yet another example.

***Example 3.3.17:*** Let $\{Q, +, \odot\}$ be a near ring. $\{Q^o, +, \odot\}$ is a seminear ring and $Q^o \subseteq Q$. Hence $\{Q, +, \odot\}$ is a definite special near ring.

We prove by examples that all near rings in general need not be a S-definite special near ring.

***Example 3.3.18:*** Let $\{Z_5, +, \odot\}$ be a near ring with $(Z_5, +)$ a group under addition modulo 5. Further $a \odot b = a$ for all $a \in Z_5$. $\{Z_5, +, \odot\}$ is a near ring. We see $\{Z_5, +, \odot\}$ is not a S-definite special near ring; for we cannot find a proper subset S in $Z_5$ such that $\{S, +, \odot\}$ is a seminear ring.

We just leave the following theorem as an exercise for the reader to prove.

**THEOREM 3.3.5:** *Every near ring need not in general be a S-definite special near ring.*

Now we proceed on to define the notion of Smarandache definite special subnear ring.

**DEFINITION 3.3.8:** *Let $\{N, +, \odot\}$ be a near ring. Let S be a proper subset of N, such that $\{S, +, \odot\}$ is a subnear ring of $\{N,$*

+, ⊙}. If {S, +, ⊙} is a S-definite special near ring then we call S to be a Smarandache definite special subnear ring (S-definite special subnear ring) of the near ring N.

We prove the following interesting theorem.

**THEOREM 3.3.6:** *Let {N, +, ⊙} be a near ring. Suppose {S, +, ⊙} be a S-definite special subnear ring of N; then N is a S-definite special near ring.*

*Proof:* Given {N, +, ⊙} = N is a near ring which contains a proper subset S ≠ ϕ and S ≠ N such that {S, +, ⊙} is a S-definite special subnear ring of N; so S contains a proper subset P, P ≠ ϕ and P ≠ S such that P is a seminear ring under the operation + and ⊙. Now P ⊆ S and S ⊆ N so P ⊆ S ⊆N and P is a seminear ring in N so N itself is a S-definite special near ring.

It is important to note that even if {N, +, ⊙} is a S-definite special near ring; N need not contain a S-definite special subnear ring.
    The task of finding an example to this effect is left as an exercise to the reader.

Now we proceed on to define more properties about S-definite special near rings.

**DEFINITION 3.3.9:** *Let {N, +, ⊙} be a near ring. I an ideal of N. We define I to be a S-definite special ideal of N if I contains a proper subset P, P ≠ ϕ and P ≠ I such that P is a seminear ring under the operations of N.*

**DEFINITION 3.3.10:** *Let {N, +, ⊙} be a S-definite special near ring. (P, +) be a group with ⊙ and let N be a near ring. Let μ: N × P → P; (P, μ) is called an N-group if for all p ∈ P and for all n, $n_1$ ∈ N; we have (n + $n_1$) p = np + $n_1$p and (n$n_1$)p = n($n_1$p).*
    *We call the N-group to be a Smarandache definite special N-group if the N-group is a S-definite special group. So a N-*

*subgroup of P will be called as a S-definite special N-subgroup if the subgroup is itself a S-definite special N-group.*

We propose some problems in this direction in the last chapter of this book.

Now we can have nontrivial class of S-definite special near rings. This is got by constructing group near rings using S-definite special near rings and any group. For more about group near rings please refer [106].
    The notion of S-definite special seminear rings is identical with S-definite special near rings introduced in this book.
    One can define S-definite special near ring homomorphisms as in [106]. We have suggested problems in this in the last chapter of this book.
    Further the interested reader can see that the S-definite special near rings can be used to construct special automatons which can do sequential operations provided the S-definite special near ring N has finite subsets in N such that (S, + ) and (S, ⊙) can be generated as free semigroups. In view of this we give the applications of S-definite special near rings in finite automaton.

Now we proceed on to define S-special near automaton using S-definite special near rings.

**DEFINITION 3.3.11:** *Let {N, +, ⊙} be a S-definite special near ring. Let P be a finite proper subset in N such that P generates a free semigroup under addition and ⟨P, +⟩ ⊆ {N, +}. Then the semiautomaton can be constructed using the set P with respect to + as the alphabets, for any set of finite states Z. Now we have $\mu: Z \times P \to Z$ is a semiautomaton; since P generates a free semigroup ⟨P⟩ and ⟨P⟩ = P is contained in N. We see $\overline{\mu}: Z \times P \to Z$ gives a semiautomaton which can work on a sequence of elements as an input alphabet.*

Now a S-definite special near ring endowed with such a semiautomaton will be known as the Smarandache near definite special semiautomaton. As we are unaware of the fact that whether near rings having only finite number of elements can ever be a S-definite special near ring.

We proceed on to define S-near definite special automaton.

**DEFINITION 3.3.12:** *Let {N, +, ⊙} be a S-definite special near ring. Let P be a proper finite subset of N. Suppose (P, +) generates a free semigroup $\overline{P_+}$ under + and {P, ⊙} generates a free semigroup $\overline{P_\odot}$ under ⊙ so that {P, +, ⊙} is a seminear ring in {N, +, ⊙} then an automaton associated with P; and $\overline{P_\odot}$ and $\overline{P_+}$ denotes the automaton for some set of states with output alphabet S again to a proper finite subset of N with $\overline{S_+}$, a free semigroup generated by S with respect to + and $\overline{S_\odot}$ a free semigroup generated by S with respect to ⊙ and $\{\overline{S}, \odot, +\}$ is a seminear ring then for any set of finite state Z we with a map $\mu$ : Z×P → P and $\lambda$ : Z×P → S and $\overline{\mu} : Z \times \overline{P_+} \to \overline{P_+}$ and $\overline{\lambda} : Z \times \overline{P_\odot} \to \overline{S_\odot}$ where $\overline{\mu}$ and $\overline{\lambda}$ are extensions of $\mu$ and $\lambda$ respectively; then we call $\{Z, \overline{P_+}, \overline{P_\odot}, \overline{S_+}, \overline{S_\odot}, \overline{\mu}, \overline{\lambda}\}$ to be the Smarandache near definite special automaton (S-near definite special automaton) associated with the S-definite special near ring N.*

*Remark:* Unless we have N to be a S-definite special near ring we cannot think of a S-near definite special automaton associated with N.

Thus one of the applications of S-definite special near ring is its use in construction of finite machines.

It has become important to mention here that associated with a S- definite special near ring N we can have several S-near definite special automatons and semiautomatons associated with it.

Thus this is the main advantage of using S-definite special near rings in constructing S-near definite special automatons. So a single algebraic structure can pave way for several finite machines. Now one may wonder how a S-definite special commutative near ring is. To this end we define the following.

**DEFINITION 3.3.13:** *Let N be a S-definite special near ring. If we have atleast one seminear ring in N under the same operations of N to be commutative, then we call N to be a S-definite special near ring or trivially S-definite special commutative near ring. Further we see all S-definite special commutative near rings need not be commutative but only S-definite special commutative near rings.*

We now proceed on to define S-definite special strongly commutative near rings.

**DEFINITION 3.3.14:** *Let {N, +, ⊙} be a near ring. If every proper subset P of N which is a seminear ring is commutative; then we call (N, +, ⊙) to be a Smarandache definite special strongly commutative near ring (S-definite special strongly commutative near ring).*

We have the following important theorem.

**THEOREM 3.3.7:** *Let (N, +, ⊙) be a commutative near ring which is a S-definite special near ring then N is a S-definite special strongly commutative near ring.*

Proof is obvious from the definition.

The reader is requested to prove that a S-definite special commutative near ring in general is not a S-definite special strong commutative near ring.
    Finally it has become pertinent to mention here that we cannot construct S-definite special semirings for we see semirings are the most generalized concepts of rings and fields. Just like we do not have the notion of Smarandache groups we cannot have the notion of Smarandache definite special semirings. However we have the notion of S-definite special near rings.

Chapter Four

# SUGGESTED PROBLEMS

In this chapter we suggest over 200 problems for the reader to solve. Some of the problems are simple and easy to solve. Some of them can be treated as conjectures. The solving of the problems will certainly make the reader not only more involved but also make the reader get insight in to the subject. Further the Smarandache definite special algebraic structures try to find a weaker structure in the stronger structure unlike the S-algebraic structures.

1. Prove there exists no finite group which is a S-special definite group.

2. Give an example of a strongly commutative S-special definite group which is not a commutative S-special definite group.

3. Does the group $\{Q \setminus \{0\}, \times\}$ have a subgroup which is not a S-special definite subgroup of $\{Q \setminus \{0\}, \times\}$.

4. Find S-special definite groups which are simple!

5. Give an example of a non commutative S-special definite group which has non trivial S-special definite normal subgroups.

6. Give an example of a S-special definite group which has only S-special definite subgroups but does not contain S-definite special normal subgroups.

7. Does the Cauchy theorems prove the non existence of finite S-special definite groups?

8. Does there exist group G which has torsion free elements yet G is not a S-special definite group?

9. Does there exist S-special definite group G which has a S-special definite minimal ideal?

10. Prove or disprove a S-special definite maximal ideal related to the semigroup of a S-definite special group is a S-special definite principal ideal!

11. Let G be a S-special definite group. Let K and H be semigroups. Find conditions for x in G or K and H so that KxH is a semigroup of G.

12. Can we have a S-special definite group G such that for semigroups H and K HxK is never a semigroup for all x ∈ G?

13. Let G be a S-special definite group H and K be semigroups of G. Does there exist x ∈ G such that H x K is a subgroup of G?

14. Characterize those S-special definite groups such that the product of every pair of semigroups is again a semigroup.

15. Illustrate by an example the product of two semigroups is not a semigroup in a S-definite special group.

16. Is in Q \ {0}, the S-definite special group the product of any two semigroups a semigroup.?

17. Is in the S-definite special group $M_{2 \times 2}$ the product of any two semigroups a semigroup?

18. Is in the S-definite special group Z under '+' the sum of any two semigroups a semigroup?

19. Can a finite field be a S-special definite field? Justify your claim.

20. Prove $\dfrac{Z_{11}[x]}{\langle x^2 + x + 1 \rangle = I} = T$ (where I is the ideal generated by $x^2 + x + 1$) is not a S-special definite field.

21. Let F be a field. Can F contain a S-special definite minimal ideal?

22. Let F be an extension field of K. Suppose F is a S-definite special field and F contains T as a subring. If G(K, F) be the group of automorphisms of K keeping F fixed and if $G_s$(T, P) (where P ⊂ F is a ring in F and T ⊂ K is a ring in K) be the group of automorphisms of T keeping P fixed. When will G(K, F) be isomorphic to ($G_s$(T, P))?

23. Let Q be a field and Q($\sqrt{2}$, $\sqrt{3}$) be the smallest field containing Q, $\sqrt{2}$ and $\sqrt{3}$; i.e. Q($\sqrt{2}$, $\sqrt{3}$) is the finite extension of Q. Let Z be the ring in Q and Z ($\sqrt{2}$, $\sqrt{3}$) be the ring generated by Z, $\sqrt{2}$ and $\sqrt{3}$. Let G (Q($\sqrt{2}$, $\sqrt{3}$), Q) be the group of automorphisms of Q ($\sqrt{2}$, $\sqrt{3}$) keeping Q fixed and $G_s$(Z($\sqrt{2}$, $\sqrt{3}$), Z) be the group of automorphisms of the ring Z ($\sqrt{2}$, $\sqrt{3}$) keeping Z fixed. Is $G_s(Z(\sqrt{2},\sqrt{3}), Z) \cong G(Q(\sqrt{2},\sqrt{3}); Q)$?

24. Determine any interesting relation between $G_s(R(a), R)$ and $G(F(a), F)$; $R \subset F$ and $R(a) \subset F(a)$ where R and R(a) are ring in F and F(a) respectively.

25. Let Q ($\sqrt{2}, \sqrt{3}, \sqrt{7}$) and Q($\sqrt{5}, \sqrt{11}, \sqrt{6}$) be two given S-definite special fields. Let Z ($\sqrt{3}$) and Z ($\sqrt{5}, \sqrt{11}$) be rings in Q ($\sqrt{2}, \sqrt{3}, \sqrt{7}$) and Q ($\sqrt{5}, \sqrt{11}, \sqrt{6}$) respectively. Give S-definite special homomorphisms $\phi_1$ and $\phi_2$ such that
   a. Ker $\phi_1$ is a trivial ideal of Z ($\sqrt{3}$).
   b. Ker $\phi_2$ is a nontrivial ideal of Z ($\sqrt{3}$).

26. For the problem 25 if Z($\sqrt{2}$) is the ring taken from Q($\sqrt{2}, \sqrt{3}, \sqrt{7}$) and Z($\sqrt{5}$) the ring taken from Q($\sqrt{5}, \sqrt{11}, \sqrt{6}$). Is it possible to find a S-definite special isomorphism from Z ($\sqrt{2}$) to Z ($\sqrt{5}$)?

27. Is $M_{2 \times 2} = \left\{ \begin{pmatrix} a & b \\ c & d \end{pmatrix} \middle| a, b, c, d \in R, ad - bc \neq 0 \right\}$ a S-special definite division ring?

28. Let K be the ring of quarternions which is a division ring. $M_{2 \times 2} = \left\{ \begin{pmatrix} a & b \\ c & d \end{pmatrix} \middle| ad - bc \neq 0; a, b, c, d \in R \right\}$ be another division ring. Construct two sets of S-special definite homomorphism from K to $M_{2 \times 2}$ so that one of them has a trivial kernel and another has a non trivial kernel.

29. Give example of a S-definite special ring.

30. Does there exist a S-definite special ring of finite order?

31. Give an example of S-definite special ring which does not contain a S-special definite ideal?

32. Give an example of a S-definite special ring R in which every S-definite special subring is a S-definite special ideal of R.

33. Give an example of a ring which is not a S-strong definite special ring!

34. Give an example of a ring other than Z which is a S-strong special definite ring.

35. Is Z[x] a S-strong special definite ring (Z[x] the polynomial ring in an indeterminate x)?

36. Is R[x] a S-strong special definite ring? (R[x] the polynomial ring with coefficients from the field of reals).

37. Does R[x] have a proper subring which is not an ideal? (We do not want R ⊆ R[x], R is a ring which is not an ideal of R[x]).

38. Give an example of a ring R in which R has no S-special definite ideals (We do not want the reader to take R to be any field it should strictly be a ring).

39. Is nZ (n a positive integer) a S-strong special definite ring? (Justify your claim).

40. Let $M_{3\times 3}$ = {M = $(m_{ij})$ | $m_{ij}$ ∈ Q} be the set of all matrices with entries from Q. Is $M_{3\times 3}$ a S-strong special definite ring?
    (a) Does $M_{3\times 3}$ contain ideals which are not S-definite special ideals?
    (b) Does $M_{3\times 3}$ contain subrings which are not S-definite special subrings?
    (c) Is every subring a S-special definite subring of $M_{3\times 3}$?

(d) Is every ideal of $M_{3\times 3}$ a S-special definite ideal of $M_{3\times 3}$?

41. Define some new and interesting properties about S-special definite rings!

42. Is $P_{2\times 2} = \left\{ \begin{pmatrix} a & b \\ c & d \end{pmatrix} \middle| a, b, c, d \in Z \right\}$ a S-strong special definite ring? Justify your claim.
    (a) Find S-definite special ideals in $P_{2\times 2}$!
    (b) Find S-definite special subrings in $P_{2\times 2}$.
    (c) Is every S-definite special subring of $P_{2\times 2}$ a S-definite special ideal of $P_{2\times 2}$?
    (d) Does $P_{2\times 2}$ contain an ideal which is not a S-special definite ideal?
    (e) Does $P_{2\times 2}$ contain a subring which is not a S-special definite subring?

43. Is $R = Z \times Z \times Z$ (the direct product of the ring of integers) = $\{(x, y, z) \mid x, y, z \in Z\}$ a S-strong special definite ring?

44. Does $R = Z \times Z \times Z = \{(x, y, z) \mid x, y, z \in Z\}$ have ideals which are not S-definite special ideals of R?

45. Does $P = Q \times Q \times Z \times R$ contain S-definite special ideals? Can P have subrings which are not S-definite special subrings?

46. Is Q[x] a S-ideally strong definite special ring? Justify your claim.

47. Can Z[x] be a S-ideally strong definite special ring?

48. Give an example of a S-ideally strong definite special ring which is not Z.

49. Is $M_{3\times 3} = \{M = (m_{ij}) \mid m_{ij} \in R$, the set of reals$\}$; a S-ideally strong definite special ring?

50. Is $P = Z \times Z \times Q$ a S-ideally strong definite special ring?

51. Can a ring $Z_n$, n an integer be a S-definite special ring?

52. Is $Z_{27}$ a S-definite special ring?

53. Is $Z_6G$(where $G = \langle g \mid g^6 = 1 \rangle$) the group ring a S-definite special ring?

54. Let ZG be the group ring of the group G where $G = D_{27} = \{a, b \mid a^2 = b^7 = 1; bab = a\}$. Find a subring in ZG which is a S-definite special ideal of ZG? Does ZG contain any subring which is not a S-definite special ideal?

55. Take the group ring $Z_{10}S_4$; the symmetric group $S_4$ over the ring of integers modulo 10. Is $Z_{10}S_4$ a S-definite special ring? Does $Z_{10}S_4$ contain a subring which is a S-definite special ring? Can $Z_{10}S_4$ have ideals which are not S-definite special ideals?

56. Is the group ring $Z_2S_3$ a S-definite special ring?

57. Can $Z_2G$ the group ring where $G = \langle g \mid g^9 = 1 \rangle$ be a S-definite special ring?

58. Does $Z_2G$ where $G = \langle g \mid g^9 = 1 \rangle$ have a subring which is a S-definite special subring?

59. Can the group ring $Z_2G$ where $G = \langle g \mid g^9 = 1 \rangle$ have a S-definite special ideal?

60. Can $Z_6G$ where $G = S_3$ have S-definite special ideals?

61. Does $Z_7S_3$ have any S-definite special ideals?

62. Does $Z_7S_7$ have any S-definite special subrings?

63. Can $Z_5 S_5$ have subrings which are not S-definite special subrings?

64. Can $Z_p S_n$, $(p, n) = 1$ have S-definite special subrings? What happens if $(p\ n) = n$ or $p$?

65. Does $ZS_5$ have any subrings which are not S-definite special subrings?

66. Is $ZS_n$ a S-strong special definite ring?

67. Is $ZS_n$ a S-ideally strong special definite ring?

68. Characterize all those subrings in $ZS_5$ which are not S-sepcial definite ideals?

69. Find all S-special definite subrings in $ZS_5$ which are not S-special definite ideals!

70. Find in ZS where S is the semigroup of all mappings of the set (1 2 3) to (1 2 3) under composition of mappings a S-definite special subring which is not a S-definite special ideal of ZS.

71. Let S(5) denote the set of all mappings of (1 2 3 4 5) to itself; under composition of mappings S(5) is a semigroup with unit. Let Z be the ring of integers Z(S(5)) the semigroup ring.

    (a) Does Z(S(5)) have a subring which is not a S-special definite ring?
    (b) Does Z(S(5)) contain an ideal which is not a S-special definite ideal?
    (c) Give an ideal in Z(S(5)) which is a S-special definite ideal.

72. Is $Z_3 S(5)$ a S-definite special ring?

73. Can $Z_5S(5)$ have S-definite special ideals?

74. Does $Z_6S(5)$ contain S-definite special subrings?

75. Can the semigroup ring $Z_nS(m)$; $(n, m) = 1$ be a S-definite special ring?

76. Is $Z_{18} S(6)$ the semigroup ring a S-definite special ring?

77. Does there exist a S-definite special ideal in $Z_6 (S(6))$?

78. Give an ideal in $Z_4S(6)$ which is not a S-definite special ideal.

79. Let $ZS(5)$ be the semigroup ring. Is every ideal in $ZS(5)$ a S-special definite ideal? If so find atleast one ideal; which is a S-special definite ideal.

80. Let $ZS(3)$ be the semigroup ring. Is every subring in $ZS(3)$ a S-special definite subring?

81. Let $ZS(4)$ be semigroup ring. Is in $ZS(4)$ every S-special definite subring a S-special definite ideal?

    [Can the problems 79, 80, 81 be studied for $ZS(n)$; $n < \infty$].

82. Let $M_{2 \times 2} = \{M = (m_{ij}) \mid m_{ij} \in Z\}$ be the ring of $2 \times 2$ matrices with entries from Z. Does $M_{2 \times 2}$ contain ideals which are not S-definite special ideals? Give a S-definite special subring in $M_{2 \times 2}$ which is not a S-definite special ideal.

83. Let $P_{3 \times 3} = \{P = (p_{ij}) \mid p_{ij} \in Q\}$ be the ring of $3 \times 3$ matrices with entries from Q. Is $P_{3 \times 3}$ a S-strong definite special ring? Justify your claim! Does $P_{3 \times 3}$ contain S-definite special subrings?

84. $Q[x]$ be the polynomial ring in the variable x. Can $Q[x]$ have S-definite special subrings which are not S-definite special ideals?

85. Is every S-definite special subring which is an ideal a S-definite special ideal in a ring?

86. Can $Z_p[x]$ a polynomial ring over $Z_p$ be a S-definite special ring?

87. Will $Z_n[x]$, n a composite number have S-definite special ideals?

88. Find in $Z[x_1, x_2, x_3]$ subrings which are not S-definite special subrings.

89. Given $R[x_1, x_2, \ldots, x_n]$ is a S-definite special ring where, R is just a commutative ring with unit. Can R be a S-definite special ring?

90. Is the class of rings $Z_n[x]$; n a composite number only be a S-ring and not a S-definite special ring?

91. What can one say about the class of rings $Zp[x_1, x_2, \ldots, x_n]$, p a prime and $Z_p$ the prime field of characteristic p and $x_1, x_2, \ldots, x_n$ are n-variables and $Z_p[x_1, x_2, \ldots, x_n]$ are polynomial rings in n-variables $x_1, x_2, \ldots, x_n$ with coefficients from $Z_p$.

92. Prove R[x] where R is the field of reals is a S-ring as well as S-definite special ring.

93. Is Z[x] a S-ring as well as a S-definite special ring?

94. Is $T_{p \times p} = \{T = (t_{ij}) \mid t_{ij} \in Z_n$ the ring of integers modulo n$\}$; a ring under matrix multiplication and matrix addition a S-ring and a S-definite special ring? Justify your claim!

95. Prove $Q[x_1, x_2, x_3]$ is a S-ring as well as a S-definite special ring. Find ideals in $Q[x_1, x_2, x_3]$ which are not S-definite special ideals.

96. Characterize those rings R which are S-rings as well as S-definite special rings.

97. Characterize those rings which are only S-rings and never a S-definite special ring.

98. Characterize those rings which are S-definite special rings but never a S-ring.

99. Does S-double strong rings have infinite cardinality?

100. Can we have S-doubly strong rings with finite number of elements in them?

101. Does there exists finite rings which are S-doubly strong rings?

102. Are all infinite rings S-doubly strong rings? Justify your claim!

103. Does there exists S-special definite rings which has finite number of elements in them?

104. Can we have any finite set which can be transformed into a S-definite special algebraic structure by defining algebraic operations on them appropriately?

105. Can a finite group i.e., a group of finite order be a S-definite special group?

106. Can a finite ring i.e., a ring with only finite number of elements be a S-definite special ring?

107. Can we have S-special definite prime fields (which are prime) to be a field which is not a prime field?

108. Prove every prime field which is a S-special definite field is a S-special definite prime field!

109. Can $Z_{19}$ be a S-special definite prime field? Justify your claim!

110. Is $Z_{23}$ a S-special definite field?

111. Can we have finite characteristic field to be a S-definite special field?

112. Suppose F is a field of characteristic p. Can F be a S-definite special field? Justify your answer.

113. Is $\dfrac{Z_{19}[x]}{\langle x^2 + x + 4 \rangle}$ a S-special definite field?

114. Is $\dfrac{Z_5[x]}{\langle x^2 + x + 1 \rangle}$ a S-special definite field?

115. Is $\dfrac{Q[x]}{\langle x^2 + 1 \rangle}$ a S-special definite field? Note $\langle x^2 + 1 \rangle$ denote the ideal generated by $x^2 + 1$ so in general $\langle p(x) \rangle$ denotes the ideal generated by the polynomial p(x).

116. Is $\dfrac{Q[x]}{\langle x^2 + 16 \rangle}$ a S-definite special prime field?

117. Is $Q(\sqrt{2})$ the field containing Q and $\sqrt{2}$ a S-definite special prime field?

118. Give examples of S-special definite fields which are not S-definite special prime fields.

119. Can you prove we have infinite collection of S-special definite fields which are not S-special definite prime fields?

120. Let $Q[x_1, x_2]$ be a vector space over Q. Is $Q[x_1, x_2]$ a S-definite special vector space over Q?

121. Can $M_{m \times n} = \{(m_{ij}) \mid m_{ij} \in Z\}$ be a S-definite special vector space over Q?

122. Is Q[x] a S-definite special vector space over R?

123. Is R[x] a S-definite special vector space over Q?

124. Can $Q \times Q \times Q \times Q = V$ be a S-definite special vector space over Q? Find a basis of V which is not a S-definite special basis of V.

125. Can every vector space V over a field F be a S-definite special vector space?

126. Does their exists a vector space defined over the field of characteristic zero which is not a S-special definite vector space?

127. Let $M_{3 \times 4} = \{(m_{ij}) = M \mid m_{ij} \in R, 1 \le i \le 3 \text{ and } 1 \le j \le 4\}$ be a vector space over R. Is $M_{3 \times 4}$ a S-definite special vector space? Give a basis for $M_{3 \times 4}$ which is a S-definite special basis.

128. Find the dimension of all the basis of all the semivector spaces over the vector space $M_{3 \times 2} = \{(m_{ij}) = M \mid 1 \le i \le 3; 1 \le j \le 2 \text{ and } m_{ij} \in Q\}$ over Q. How many semivector spaces are in $M_{3 \times 2}$? Is $M_{3 \times 3}$ a S-definite special vector space of finite dimension?

129. What is the S-definite special dimension of the vector space where the vector space $Q \times Q \times Q \times Q \times Q = V$ over Q?

130. Find the S-definite special basis for $M_{3 \times 3} = \{(m_{ij}) = m \mid m_{ij} \in Q; 1 \le i, j \le 3\}$ over Q.

131. Can any semivector space obtained from the above problem have infinite number of basis?

132. Let $V = Q \times Q \times Q \times Q \times Q$ be a vector space over Q. Is V a S-definite special vector space? How many semivector spaces are in V? Find the number of S-definite special basis? Is the basis of the semivector spaces in V unique?

133. What is the structure of $L_p(V, W)$? Here V and W are S-definite special vector spaces over a field F and $L_p(V, W)$ denotes the collection of all S-definite special linear transformations of S-definite special vector spaces.

134. Suppose $L_p(V^s, W^{S_1})$ denotes the set of all S-definite special linear transformations of the S-definite special vector spaces over F restricted to the set of all linear transformations from the semivector space S to $S_1$, what is the algebraic structure enjoyed by $L_p(V^s, W^{S_1})$, under the sum of the transformations and composition of transformation when ever defined?

135. Let $V = M_{3 \times 2} = \{M = (m_{ij}) \mid m_{ij} \in Q; 1 \leq i \leq 3 \text{ and } 1 \leq j \leq 2\}$ be a S-definite special vector space over Q and $W = M_{2 \times 2} = \{M = (m_{ij}) \mid m_{ij} \in R, 1 \leq i, j \leq 2\}$ be a S-definite special vector space over Q. Give 2 nontrivial S-special definite linear transformations from V to W. Can we ever have a S-definite special transformation which is 1-1?

136. Let $V = R \times R$ over the field Q and $W = Q \times Q \times Q$ over the field Q be two S-definite special linear vector spaces. For $S = Z^\circ \times Z^\circ$ over $Z^\circ$ and $S_1 = Z^\circ \times Z^\circ \times Z^\circ$ over $Z^\circ$, the S-special transformation given by $T(x, y) = (x + y, 3x, x - y)$. Find Ker T and the matrix of S-special definite linear transformation related with T.

137. Let $V = R \times R \times R \times R$ be a S-definite special vector space over Q. Find two proper S-special definite subspaces of V.

138. Does a S-definite special vector space contain two S-special definite subspaces $W_1, W_2$ such that $W_1 \cap W_2 = \{0\}$?

139. Obtain any other interesting results about S-definite special subspaces of a vector space.

140. Can we define for any S-definite special vector space V a S-special definite subspaces W and $W^\perp$ such that $W + W^\perp = V$?

141. Bring out the properties enjoyed by the S-special definite inner product on a vector space V. Does it depend on the semivector space choosen in V? Illustrate this explicitly by examples.

142. Show that the S-speduo special standard inner product on V $= Q \times Q \times Q \times Q$ over Q on the S-definite special vector space V for which the semivector space is taken as $W = Z^\circ \times Z^\circ \times Z^\circ \times Z^\circ \times Z^\circ$ over $Z^\circ$, one cannot normalize each of the S-special definite vectors as unit vectors. Prove $\{0\} = (0\ 0\ 0\ 0)$ is the only vector which is orthogonal under the S-definite special inner product.

143. Let $V = Q \times Q \times Q \times Q \times Q$ be a S-definite special vector space over Q. Let $W = Z^\circ \times Z^\circ \times Z^\circ \times Z^\circ \times Z^\circ$ be a semivector space of V over the semifield $Z^\circ$. Suppose $(\alpha / \beta) = 3x_1y_1 - x_2y_2 + x_3y_3 - x_4y_4 - x_5y_5$ be an inner product on the vector space V where $\alpha = (x_1, x_2, x_3, x_4, x_5)$ and $\beta = (y_1, y_2, y_3, y_4, y_5)$. Can the restriction of this inner product be a S-special definite inner product (relative to the semivector space $W \subseteq V$) of the vector space V.

144. Let G be a S-definite special group and K any S-definite special field. Is the group ring a S-definite special ring?

145. Give an example of a group ring which is not a definite special ring.

146. Let G be a S-definite special group and K any field or a commutative ring with unit. Is the group ring KG a S-definite special ring?

147. Let G be a S-semigroup and K any field. Can the S-semigroup ring KG be a S-ring?

148. Let R be a S-definite special ring with 1 and G any group. Is the group ring RG a S-definite special ring?

149. Does there exists any S-definite special linear algebra over a field of finite characteristic p ($p < \infty$, $p > 0$)?

150. Is $Z_7[x]$ the linear algebra over $Z_7$ a S-definite special linear algebra?

151. Let $M_{3\times 3} = \{M = (m_{ij}) | m_{ij} \in Z_5; 1 \leq i, j \leq 3\}$ be the collection of all $3 \times 3$ matrices with entries from $Z_5$. Is $M_{3\times 3}$ a S-definite special linear algebra over $Z_5$? Justify your claim!

152. Let $V = Q \times Q \times Q \times Q \times Q$ be a linear algebra over Q Find all the semilinear algebras contained in V. Find the S-definite special basis of V over Q.

153. Let $V = R \times R \times R \times R$ be a linear algebra over Q. Is V a S-special definite linear algebra over Q? Is V a finite dimensional S-special definite linear algebra?

154. Let $V = R \times R \times R \times R$ be a linear algebra over R and $W = Q \times Q \times Q$ be a linear algebra over Q; T: V → W be a map such that T(x y z w) = (x y z) where P = R° × R° × R° × R° is a semilinear algebra over R° and S = Q° Q° Q° is a semilinear algebra over Q°. Can T be a S-definite linear algebra transformation for the semilinear algebras P and S? Substantiate your answer! Is T defined from V to W a S-definite linear algebra transformation? Can T be a linear algebra transformation?

155. Suppose in the above problem if V is also defined over Q. Is T a linear algebra transformation? Is T S-definite linear algebra transformation?

156. Let $V = Q \times Q \times Q \times Q$ be a S-definite special vector space over Q. How many S-definite special converging transformation on V can be defined?

157. Let $M_{n \times n} = \{M = (m_{ij}) \mid m_{ij} \in R\}$ be a S-definite special vector space over Q. Define a S-definite special converging transformation on $M_{n \times n}$.

158. Let $V = Q[x]$ be a S-definite special vector space over Q. Does there exists a transformation from V to V whose restriction can be S-definite special converging transformation on V. How many such transformations can be defined?

159. Let V be a S-definite special vector space over a field F. Is every S-definite special converging transformation be extended to be a linear transformation?

160. Let $V = M_{n \times n} = \{M = (m_{ij}) \mid m_{ij} \in Q; 1 \le i \le m; 1 \le j \le n\}$ be a S-definite special vector space over Q. Let $W = M^s_{m \times n} = \{M = (m_{ij}) \mid m_{ij} \in Z^\circ, 1 \le i \le m, 1 \le j \le n\}$ be a semivector space over $Z^\circ$. Define a S-definite special diverging transformation $T_D$ on V. Define a S-definite special converging $T_c$ transformation on V. Can ever $T_D$ and. $T_c$ be related in any way?

161. Obtain some interesting properties about $T_D$ and $T_c$.

162. Define $T_D$ which is such that $T_D^{-1}$ is a S-definite special linear transformation!

163. Is it possible to define $T_c$ such that $T_c^{-1}$ exists and T is a S-definite special linear transformation?

164. Given $V = Q \times Q \times Q \times Q$ is the S-definite special linear transformation. Define $T_D$ and $T_c$ such that $T_D^{-1}$ exists? Is it possible to define $T_c^{-1}$ on V?

165. Obtain some interesting results on V when V is both S-definite special vector space as well as S-vector space.

166. Give an example of a S-definite special near ring, which has no S-definite special subnear ring.

167. Give an example of a S-definite special near ring which has a subnear ring which is not a S-definite special subnear ring.

168. Can we have S-definite special near rings which has only finite number of elements in them?

169. Is $\{Z_{10}, +, \odot\}$ with a $a \odot b = a$, a S-definite special near ring?

170. Can $\{Z_{17}, +, \odot\}$ be a S-definite special near ring?

171. Give an example of a S-definite special ideal in a S-definite special near ring.

172. Is every S-definite special ideal, an ideal of the S-definite special near ring?

173. Can every ideal of a near ring be a S-definite special ideal of a near ring?

174. Can $\{Z, +, \odot\}$ have S-definite special ideals?

175. Can $\{Z_{11}, +, \odot\}$ have S-definite special ideals?

176. Does $\{Z_{18}, +, \odot\}$ have S-definite special ideal?

177. Is every ideal of $\{Z, +, \odot\}$ a S-definite special ideal of $\{Z, +, \odot\}$?

178. Can the near ring $\{Z, +, \odot\}$ have S-definite special N-group associated with it?

179. Can one say every N-group of any near ring will be a S-definite special N-group?

180. Is it true that every N-group of a S-definite special near ring of a S-definite special N-group?

181. Can a N-group associated with a S-special definite near ring be not a S-definite special N-group?

182. Find any relation or interesting conditions/ relation, between S-definite special near ring N and the S-definite special N-group of N.

183. Can one prove if $\{I_k\}$; $k \in K$ denotes the collection of all S-definite special ideals of a near ring N then the following two condition are equivalent
    a. The S-definite special ideal of N generated by $\bigcup_{k \in K} I_k$
    b. The set of all finite sums of the elements of the $I_k$'s.

184. Can the set of n × n matrices with entries from Q be made into a near matrix ring which is a S-definite special ring?

185. Define S-definite special right ideal of a near ring and illustrate it with examples.

186. Define S-definite special minimal right ideal.

187. Can the near ring $\{Z, +, \odot\}$ have S-definite special minimal right ideal?

188. Define S-definite special equiprime left ideal of a near ring N.

189. Obtain some interesting results on S-definite special near rings.

190. Illustrate the above definition by an example.

191. Define the notion of S-definite special n ideal near ring.

192. Give an example of a S-definite special n-ideal near ring.

193. Can $\{Z, +, \odot\}$ be a S-definite special n-ideal near ring?

194. Prove if the near ring N is a S-definite special near ring then for any group G the group near ring NG is a S-definite special near ring.

195. Suppose N is not a S-definite special near ring can for any group G the group near ring NG be a S-definite special near ring?

196. Suppose we take G to be a S-definite special group and N to be any near ring not in particular a S-definite special near ring. Can the group near ring NG be a S-definite special near ring?

197. Give examples of S-definite special near rings which are group near rings.

198. Give examples of group near rings which are not S-definite special near rings.

199. Can the group near ring $Z_{12}G$ where $G = \langle g \mid g^6 = 1 \rangle$ be a S-definite special near ring? (Here $Z_{12}$ is endowed with usual modulo 12 addition and multiplication is $a \odot b = a$ for all $a \in Z_{12}$).

200. Define the notion of S-definite special equiprime near ring. Will all S-definite special near rings be S-definite special equiprime near rings? Justify your claim!

201. Obtain conditions under which S-infra near ring and S-definite special near ring will be identical! Is it possible?

202. Can a S-definite special near ring be S-infra near ring?

203. Let Z be the set of integers. Let $\{N = Z \times Z \times Z, +, \odot\}$ be a near ring. Using the set $S = \{(1, 2, 3)\} \cup (0\ 0\ 0)$ generate a S-near definite special automation for any arbitrary set of states

$$\tilde{Z} = \{Z_1, Z_2, \ldots, Z_n\} \mid n < \infty\}.$$
$$P = \{(4, 7, 5)\} \cup (1\ 1\ 1)$$

where
$$\overline{\lambda} : Z \times \overline{P} \to \overline{P}$$
and
$$\overline{\delta} : Z \times \overline{P} \to \overline{S}$$

where
$$\overline{P} = \{\overline{P_+}, \overline{P_\odot}\}; \{\{(4Z°, 7Z°, 5Z°)\} = \overline{P_+}$$

where
$$4Z° = \{0, 4, 8, 12, 16, \ldots\},$$
$$7Z° = \{0, 7, 14, 21, \ldots\} \text{ and}$$
$$5Z° = \{0, 5, 10, 15, \ldots\}$$

and $(4Z°, 7Z°, 5Z°) = \{(x, y, z) \mid x \in 4Z°, y \in 7Z°, z \in 5Z°\}$;

$$\overline{P_\odot} = \{\overline{(4Z^+)}, \overline{(7Z^+)}, \overline{(5Z^+)}\}$$

where
$$\overline{(4Z^+)} = \{4, 4^2, 4^3, \ldots\};$$
$$\overline{(7Z^+)} = \{7, 7^2, 7^3, \ldots\} \text{ and}$$
$$\overline{(5Z^+)} = \{5, 5^2, 5^3, \ldots\}.$$

Now $\overline{S} = \{\overline{S_+}, \overline{S_\odot}\}$; worked as in case of $\overline{P}$. We can now define $\overline{\lambda}_p : \tilde{Z} \times \overline{S} \to \overline{S}$ depending on what is taken as input symbols $\overline{P}$ or $\overline{S}$.
So that

$$\bar{\mu}_p : \tilde{Z} \times \bar{P} \to \bar{S}$$

and

$$\bar{\mu}_s = \tilde{Z} \times \bar{S} \to \bar{P}.$$

Hence prove or disprove

$$\bar{A}_p = \{\bar{Z}, \bar{P}, \bar{S}, \bar{\lambda}_p, \bar{\mu}_p\}$$

will be a S-near definite special automation associated with N or

$$\bar{A}_s = \{\bar{Z}, \bar{S}, \bar{P}, \bar{\lambda}_s, \bar{\mu}_s\}$$

will be a S-near definite special automation.
Further prove or disprove

$$\bar{S}_p = \{\bar{Z}, \bar{P}, \bar{\lambda}_p\}$$

is a S-near definite special semiautomation associated with N and

$$\bar{S}_s = \{\bar{Z}, \bar{S}, \bar{\lambda}_s\}$$

a S-near definite special semiautomation associated with N.

204. Give examples of S-near definite special automation associated with any S-definite special near ring.

205. Prove S-near definite special automation associated with a near ring N exists if and only if N is a S-definite special near ring.

206. Does there exists a S-near definite special semiautomation with a near ring, which is not a S-definite special near ring?

207. Suppose (N, +, ⊙) is a near ring and P is a finite subset associated with a S-near definite special semiautomation of N will ($\bar{P}_+$, +, ⊙) be a semiring is general?

208. Is it always be possible to make a S near definite special semiautomation into a S-near definite special automation?

209. Define S-definite special τ-near ring and illustrate it by an example.

210. Can $(Z, +, \odot)$ be made into a S-definite special $\tau$-near ring for an appropriate $\tau$?

211. Give an example of a S-definite special strongly commutative near ring.

212. Give an example of a S-definite special commutative near ring.

213. Show by an example a S-definite special commutative near ring is not in general a S-definite special strongly commutative near ring.

214. Give an example of a S-definite special strongly commutative near ring, which is not a commutative near ring.

215. Can we define a special operation on $(Z, + \odot_s)$ so that $(Z, +, \odot_s)$ is only a S-definite special commutative near ring which is not a commutative near ring.

216. Can $(Z, +, \odot)$ be made into a S-definite special strongly commutative near ring which is not a commutative near ring?

217. Does there exists finite S-definite special commutative near ring which is not a commutative near ring?

218. Does there exists finite S-special definite strongly commutative near ring which is not a commutative near ring?

219. Let $M_{n \times n} = \{M = (a_{ij}) \mid a_{ij} \in Z\}$ be matrix near ring. Is this a S-definite special commutative ring?

220. Give a S-definite special near ring N which has at least six S-near definite special semiautomation associated with it.

221. Give a S-definite special near ring N which has only one S-near definite special automation associated with it.

222. Does there exist a S-definite special near ring N which has an infinite number of S-near definite special automations associated with it?

# FURTHER READING

# INDEX

# ABOUT THE AUTHOR

**Dr.W.B.Vasantha Kandasamy** is an Associate Professor in the Department of Mathematics, Indian Institute of Technology Madras, Chennai. In the past decade she has guided 12 Ph.D. scholars in the different fields of non-associative algebras, algebraic coding theory, transportation theory, fuzzy groups, and applications of fuzzy theory of the problems faced in chemical industries and cement industries.

She has to her credit 646 research papers. She has guided over 68 M.Sc. and M.Tech. projects. She has worked in collaboration projects with the Indian Space Research Organization and with the Tamil Nadu State AIDS Control Society. This is her 40$^{th}$ book.

On India's 60th Independence Day, Dr.Vasantha was conferred the Kalpana Chawla Award for Courage and Daring Enterprise by the State Government of Tamil Nadu in recognition of her sustained fight for social justice in the Indian Institute of Technology (IIT) Madras and for her contribution to mathematics. (The award, instituted in the memory of Indian-American astronaut Kalpana Chawla who died aboard Space Shuttle Columbia). The award carried a cash prize of five lakh rupees (the highest prize-money for any Indian award) and a gold medal.
She can be contacted at vasanthakandasamy@gmail.com
You can visit her on the web at: http://mat.iitm.ac.in/~wbv